\documentclass[a4paper,11pt]{article}
\usepackage[top=2.5cm,bottom=2.5cm,left=2.2cm,right=2.2cm]{geometry}
\usepackage{amsfonts}
\usepackage{mathrsfs,amscd,amssymb,amsthm,amsmath,bm,graphicx,psfrag,subfigure,url,mathtools}
\usepackage{pict2e}
\usepackage{stfloats}
\usepackage{psfrag,amsmath}
\usepackage{tikz}
\usepackage{indentfirst}
\usepackage{hyperref}
\usepackage{bookmark}
\usepackage{enumerate}
\usepackage{latexsym,euscript,epic,eepic,color}
\usepackage{multirow}
\usepackage{multicol}
\usepackage{longtable}
\usepackage{adjustbox}
\usepackage[all]{xy}
\usepackage{setspace}
\usepackage{epstopdf}
\allowdisplaybreaks
\usepackage{authblk}

\makeatletter

\renewcommand{\@seccntformat}[1]{{\csname the#1\endcsname}{\normalsize .}\hspace{.5em}}
\makeatother

\usepackage{ifpdf}
\usepackage{indentfirst}

\def \[{\begin{equation}}
\def \]{\end{equation}}

\newtheorem{thm}{Theorem}[section]

\newtheorem{claim}{Claim}
\newtheorem{lemma}[thm]{Lemma} 

\newtheorem{problem}{Problem} 

\newtheorem{fact}{Fact}

\newenvironment{wst}
{\setlength{\leftmargini}{1.5\parindent}
 \begin{itemize}
 \setlength{\itemsep}{-1.1mm}}
{\end{itemize}}
\begin{document}
\baselineskip=0.23in
\title{\bf Connected triangle-free planar graphs whose second largest eigenvalue is at most 1\thanks{S.L. financially supported by the National Natural Science Foundation of China (Grant Nos. 12171190, 11671164) and the Special Fund for Basic Scientiﬁc 
Research of Central Colleges (Grant No. 30106240144).\\[1pt] \hspace*{5mm}{\it Email addresses}: chengkunmath@163.com (K. Cheng),\  li@ccnu.edu.cn (S. Li)}}
\author[1]{Kun Cheng}
\author[,1,2]{Shuchao Li\thanks{Responding author}}
\affil[1]{Faculty of Mathematics and Statistics, and Hubei Key Lab--Math. Sci.,\linebreak Central China Normal University, Wuhan 430079, China}
\affil[2]{Key Laboratory of Nonlinear Analysis \& Applications (Ministry of Education),\linebreak Central China Normal University, Wuhan 430079, China}
\date{\today}
\maketitle
\begin{abstract}
Let $\lambda_2$  be the second largest eigenvalue of the adjacency matrix of a connected graph. In 2023, Li and Sun \cite{LiSun1} determined all the connected $\{K_{2,3}, K_4\}$-minor free graphs whose second largest eigenvalue $\lambda_2\le 1$. As a continuance of it, in this paper we completely identify all the connected $\{K_5,K_{3,3}\}$-minor free graphs without $C_3$ whose second largest eigenvalue does not exceed 1. 
This partially solves an open problem posed by Li and Sun \cite{LiSun1}: Characterize all connected planar graphs whose second largest eigenvalue is at most $1.$ Our main tools include the spectral theory and the local structure characterization of the planar graph with respect to its girth.
\vspace{0.2cm}

\end{abstract}

\section{\normalsize Introduction}
In this paper, we consider simple and undirected graphs. For graph theoretic notation and terminology not defined here, we refer to \cite{bondy1,Brouwer1}. 

The study of the second largest eigenvalue, $\lambda_2$, of a graph has been shown to be very significant, and has many applications in measuring the pseudo-randomness of the graph (see \cite[Chapeter 9]{Alon}). 
One important direction on the second largest eigenvalue of a graph is to study its multiplicity. 
Jiang, Tidor, Yao, Zhang and Zhao~\cite{Jiang} determined an upper bound on the multiplicity of $\lambda_2$, and they applied it, as a key tool, to solve an old geometry problem on equiangular lines. 
As a continuance of it, 
Chen and Hao~\cite{CGT} established an upper bound on the multiplicity of the second largest eigenvalue of a planar graph and a outplanar graph, respectively.

Another important direction on the study of the second largest eigenvalue is to determine the graphs whose second largest eigenvalue does not exceed a comparatively small constant. The so-called Hoffman problem (see \cite{Cve1982}) asks for a characterization of connected graphs whose second largest eigenvalue is at most 1. This problem attracts many researchers' attention in the past four decades \cite{Stanic4}. It is a very difficult problem, and not fully resolved until now. 

To the best of our knowledge, the first paper that concerns graphs with $\lambda_2\leq 1$ was published in 1982 by Cvetkovi\'{c}~\cite{Cve1982}.
Neumaier~\cite{Neu82} gave an algorithm to decide whether $\lambda_2$ of a tree exceeds an arbitrary bound or not.
Petrovi\'{c}~\cite{Pe1} determined all bipartite graphs with $\lambda_2\leq 1$.
Then, Xu and Shao~\cite{XuShao2} identified all non-bipartite graphs with girth at least $4$ satisfying $\lambda_2< 1$.
Xu~\cite{XuGH1} characterized all unicyclic graphs with $\lambda_2\leq 1$.
Guo~\cite{GuoSG} identified all bicyclic graphs with $\lambda_2\leq 1$.
Li and Yang~\cite{Li2} characterized all tricyclic graphs satisfying $\lambda_2\leq 1$.
Petrovi\'{c} and Mileki\'{c}~\cite{Pe2,Pe3} characterized all line graphs and generalized line graphs with $\lambda_2\leq 1$.
Gao and Huang~\cite{HQX} determined all the generalized $\theta$-graphs with $\lambda_2\leq 1$.
Stani\'{c}~\cite{Stanic2,Stanic3} characterized all nested split graphs and regular graphs with $\lambda_2\leq 1$.
Cheng, Gavrilyuk, Greaves and Koolen~\cite{Cheng1} identified the connected biregular graphs with precisely three distinct eigenvalues and $\lambda_2\leq 1$.
Later, Cheng, Greaves and Koolen~\cite{ChengXM} completely characterized all connected graphs with exactly three distinct eigenvalues satisfying $\lambda_2\leq 1$. 
Note that a graph is called $\mathcal{H}$-free if it does not contain any element in $\mathcal{H}$ as an induced subgraph.
In 2021, Liu, Chen and Stani\'{c}~\cite{Liu21} characterized all the connected 
$\{K_{1,3}, K_5- e\}$-free graphs with $\lambda_2\leq 1$. They also characterized all the diamond-free graphs with $\lambda_2\leq 1$; see \cite{Liu24}. 
For more results on this topic, we refer to \cite{CaoHong,GGuo,JCTA,ShuJL1,Stanic1,Zhang1}.

There are many results characterizing $\mathcal{H}$-free graphs with $\lambda_2 \le 1.$ 
It is interesting to develop the research on doing minor theory together with $\lambda_2.$ Motivated by \cite{CGT,Jiang,Liu21}, Li and Sun \cite{LiSun1} completely identified all the connected $\{K_{2,3}, K_4\}$-minor free graphs whose second largest eigenvalue does not exceed 1. That is, they characterized all the connected outerplanar graphs satisfying $\lambda_2\leq 1$. Furthermore, all the maximal outerplanar graphs with $\lambda_2\leq 1$ are also identified in \cite{LiSun1}. Li and Sun~\cite{LiSun1} proposed the following interesting problem.
\begin{problem}\label{Pro1}
Characterize all connected planar graphs whose second largest eigenvalue is at most $1.$
\end{problem}

In this paper, we focus on this open problem and determine all the connected triangle-free planar graphs with $\lambda_2\leq 1$. 

\noindent{\bf Outline of the paper}\ \ The remaining sections are organized as follows:
In Section~\ref{s2}, we give some necessary preliminaries. 
In Section~\ref{s3}, we show our main result. Some concluding remarks are given in Section~\ref{s4}. 

\section{\normalsize Preliminaries}\label{s2}
In this section, we give some necessary definitions and useful lemmas. 

Let $G=(V(G), E(G))$ be a simple graph with vertex set $\{v_1, \ldots, v_n\}$ and edge set $E(G)$. Then $n(G):=|V(G)|$ and $e(G):=|E(G)|$ are called the {\it order} and {\it size} of $G$, respectively. The \textit{adjacency matrix} of $G$ is an $n\times n$\ 0-1 matrix $A(G)=[a_{i j}]$ with $a_{ij}=1$ if and only if $v_i$ and $v_j$ are adjacent. As $A(G)$ is real symmetric, its eigenvalues are real and can be arranged as $\lambda_1(G)\geq \lambda_2(G)\geq \cdots \geq \lambda_n(G)$.

We call a graph $G$ {\em empty} if $E(G)=\emptyset$. For $u, v\in V(G)$, we denote by $u\sim v$ if $u$ and $v$ are adjacent, and $u\not\sim v$ otherwise. For $v\in V(G)$, let $N_G(v)$ and $d_G(v)$ be the neighborhood and the degree of $v$ in $G$, respectively.
For $S\subseteq V(G)$, we denote by $G[S]$ the {\em subgraph} of $G$ {\em induced by} $S$. 
Let $H\subset G$ denote that $H$ is an induced subgraph of a graph $G$.
We write $N_S(v)$ for $N_G(v)\cap S$. Then $d_S(v):=|N_S(v)|$ is the {\em degree} of $v$ {\em in} $S$. A vertex $v$ is called a {\em pendant vertex} of $G$ if $d_G(v)=1$. 
Let $X,Y\subseteq V(G)$ be two disjoint vertex sets, then the set of edges with one end in $X$ and the other in $Y$ is denoted by $E(X,Y)$. We use $e(X,Y)$ to denote $|E(X,Y)|.$

We write $P_n$, $C_n$ and $K_n$ for the path, the cycle and the complete graph on $n$ vertices, respectively. Denote by $K_{s,t}$ the complete bipartite graph with $s$ vertices in one colour class and $t$ vertices in the other. If $G$ contains at least one cycle, then the {\em girth} of $G$, denoted by $g(G)$, is the length of its shortest cycle. 

\begin{figure}[!ht]
\begin{center}
\includegraphics[scale=0.9]{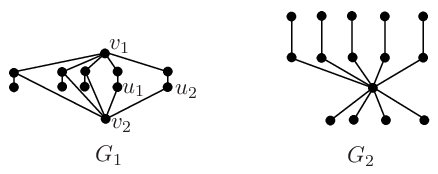}
\caption{Graphs $G_1$ and $G_2$.}\label{Fig1}
\end{center}
\end{figure}

Let $\mathcal{F}_1$ be the class of planar graphs obtained from $K_{2,s}$ ($s\geq 3$) in the following way: There exist two vertices, say $v_1, v_2$, of degree $s$  and $s$ vertices, say $u_1, u_2, \ldots, u_s$, of degree 2 in $K_{2,s}$.
We first subdivide $v_1u_1$ and $v_1 u_2$,
then attach a pendant edge to each vertex of $\{u_3,\ldots, u_s\}$ (see graph $G_1$ in Figure~\ref{Fig1}). Clearly, $G_1$ is a graph of $\mathcal{F}_1$.

Let $\mathcal{F}_2$ be the class of trees obtained from $K_{1,t}$ in the following way: Choose $k$ ($k\leq t$) pendant vertices of $K_{1,t}$. Then attach a pendant edge to each of them (see graph $G_2$ in Figure~\ref{Fig1}). It is easy to see that $G_2$ is in $\mathcal{F}_2$. 

One easily see that, for any graphs $G\in \mathcal{F}_2$, there exists a graph $G'\in \mathcal{F}_1$, such that $G\subset G'$. So, we use $\mathcal{F}_2\subset \mathcal{F}_1$ to denote this relationship. 

An {\em independent set} in $G$ is a set of vertices no two of which are adjacent. 
Let $G$ be a graph, and let $\mathcal{H}$ be a set of graphs, we say that $G$ is {\em $\mathcal{H}$-free} if it does not contain any element in $\mathcal{H}$ as an induced subgraph. A {\em minor} of $G$ is a graph obtainable from $G$ by means of a sequence of vertex deletions, edge deletions and edge contractions. A graph is said {\em $\mathcal{H}$-minor free}, if it does not contain any element in $\mathcal{H}$ as a minor. In particular, for a graph $H$, we usually write $H$-free (resp., $H$-minor free) instead of $\{H\}$-free (resp., $\{H\}$-minor free). 

A graph is {\em planar} if it can be drawn in the plane in such a way that edges meet only at points corresponding to their common ends.
A classical result in graph theory, known as the Kuratowski's Theorem, states that
a graph is planar if and only if it is $\{K_5,K_{3,3}\}$-minor free (see \cite[Theorem~10.32]{bondy1}).

\begin{figure}[!ht]
\begin{center}
\includegraphics[scale=1.1]{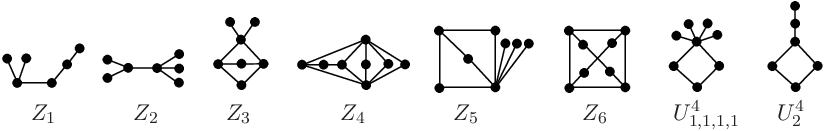}
\caption{The forbidden graphs.}\label{Fig2}
\end{center}
\end{figure} 

In 1982 Neumaier~\cite{Neu82} determined all trees with $\lambda_2\le 1$.
\begin{lemma}[\cite{Neu82}]\label{lemT}%
For a tree $T$, $\lambda_2 (T)\leq 1$ holds if and only if $T$ is an induced subgraph of some graph belonging to
$\mathcal{F}_2\cup G_0,$ where $G_0$ is obtained by attaching two pendant edges to a pendant vertex of $K_{1,3}$.\footnote{Note $\mathcal{F}_2\subset \mathcal{F}_1$ and $G_0\subset \mathcal{F}_1$.}
\end{lemma}

Let $Z_1,\ldots,Z_6,U^4_{1,1,1,1}, U^4_2$ be the graphs illustrated in Figure~\ref{Fig2}.
By a direct computation we get the following lemma.
\begin{lemma}\label{lem31}
If $G\in \{Z_1,\ldots,Z_6,U^4_{1,1,1,1}, U^4_2,P_6\}$, then $\lambda_2(G)>1$.
\end{lemma}

Recall that the definition of $\mathcal{F}_1$ depends on the parameter $s$  with $s\geq 3$.
\begin{lemma}\label{lemF1}
If $G\in \mathcal{F}_1$, then $\lambda_2(G)= 1$.
\end{lemma}
\begin{proof}
If $G\in \mathcal{F}_1$, then its spectrum consists of one copy of $-2$, $s-2$ copies of $-1$, $s$ copies of $1$ and the roots of $f(x)=x^3-2sx+x-2$. Since $f(-\infty)<0$, $f(-1)=2s+4$, $f(1)=-2s<0$ and $f(+\infty)>0$, we have $\lambda_2(G)=1$. 
\end{proof}
Throughout the remainder of the paper, we denote
\begin{align*}
\text{
$\mathcal{G}$:=\{$G$: $G$ is a connected triangle-free planar graph with $\lambda_2 (G)\leq 1$\}.
}
\end{align*}

The following notations are frequently used in the forthcoming discussions. Let $G$ be in $\mathcal{G}$ and let $C_g$ be the shortest cycle of $G$. Denote $T:=V(G)\setminus V(C_g)$.
For $S\subseteq V(C_g)$, $T_S:=\{v: v\in T$ with $N_G(v)\cap V(C_g) =S\}$ and $T_i$ denotes the set of vertices of $T$ which are adjacent to exactly $i$ vertices of $V(C_g)$. By these definitions, one has $T_1=\bigcup_{x\in V(C_g)} T_x$ and $T_0$ consists of the vertices non-adjacent to any vertex of $V(C_g)$. 

\begin{figure}[!ht]
\begin{center}
\includegraphics[width=0.9\textwidth]{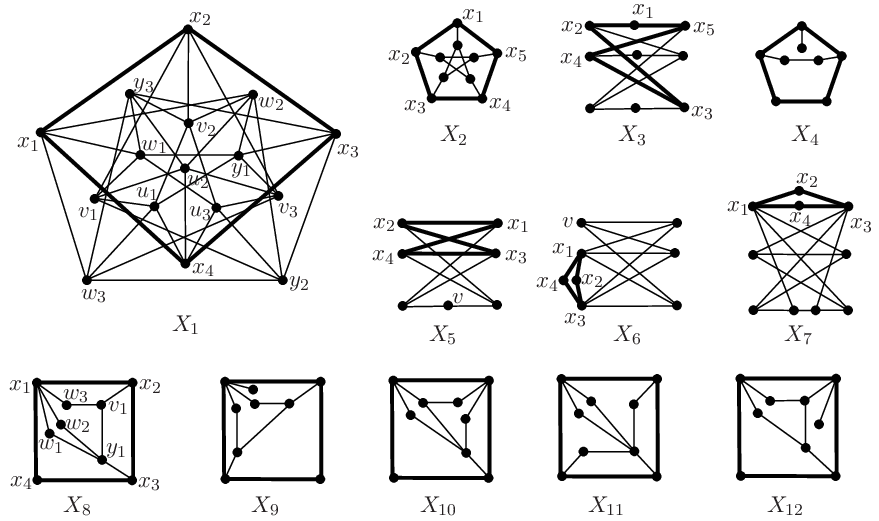}
\caption{Graphs $X_1,\ldots, X_{12}$.}\label{Fig6}
\end{center}
\end{figure}

In our context, we say that $C_g$ is a {\em base graph} of $G$. For ease of reading, we use bold lines to specify the edges of base graph.
For example, the base graph of $X_1$ (see Figure~\ref{Fig6}) is the cycle $C_4$ with the vertex set $\{x_1,x_2,x_3,x_4\}$. $X_1$ is called the {\em Clebsch graph}, a well-known unique strongly regular graph with parameters $(16,5,0,2)$ (see \cite[Theorem~10.6.4]{Godsil}).

Cauchy Interlacing Theorem (see for example \cite{Brouwer1,Godsil}) tells us that $\lambda_2\leq 1$ is a hereditary property, that is, if it holds for a graph $G$ then it holds for every induced subgraph of $G$. By Lemma~\ref{lem31}, each graph in $\mathcal{G}$ is $\{Z_1,\ldots,Z_6,U^4_{1,1,1,1}, U^4_2,P_6\}$-free. 

The following lemma will be needed.
\begin{lemma}\label{c2}
Let $G\in \mathcal{G}$ with girth $g(G)=4$, and let $C_4=x_1x_2x_3x_4x_1$ be its base graph (in what follows, the subscripts are computed modulo $4$). Then the following statements hold:
\begin{wst}
\item[{\rm (i)}] $T=T_0\cup T_1\cup T_2$, $T_2=T_{\{x_1, x_3\}}\cup T_{\{x_2, x_4\}}$, $E(T_{\{x_i, x_{i+2}\}}, T_{x_i})=\emptyset$, $G[T_{x_i}]$ and $G[T_{\{x_i, x_{i+2}\}}]$ are empty graphs, and $|T_{x_i}|\leq 3$, for $1\leq i \leq 4;$
\item[{\rm (ii)}] $T_2$ is an independent set of $G;$ 
\item[{\rm (iii)}] Every vertex of $T_2$ is adjacent to at most one vertex of $T_0;$
\item[{\rm (iv)}] Every vertex of $T_0$ is non-adjacent to any vertex of $T_1;$
\item[{\rm (v)}] If $v\in T_0$, then $d_{T_2}(v)\geq 1$. Moreover, if $d_{T_{\{x_i, x_{i+2}\}}}(v)\geq 1$, then $d_{T_2}(v)=d_{T_{\{x_i, x_{i+2}\}}}(v)$\linebreak and $d_{T_{\{x_i, x_{i+2}\}}}(v)\leq 2$. 
    Particularly, if $d_{T_{\{x_i, x_{i+2}\}}}(v)= 2,$ then $|T_2|=|T_{\{x_i, x_{i+2}\}}|=2;$
\item[{\rm (vi)}] $T_0$ is an independent set of $G$, 
hence, the neighbours of every vertex of $T_0$ belong to $T_2;$
\item[{\rm (vii)}] For $i \in\{1,2,3,4\}$, if $u\in T_{x_i}$, then $|T_{x_{i+2}}|-1\leq d_{T_{x_{i+2}}}(u)\leq 2$ and $|T_{x_j}|-2\leq d_{T_{x_j}}(u)\leq 1,$ where $j\in\{i-1,i+1\}.$
\end{wst}
\end{lemma}

\begin{proof}
(i)\ Since $g(G)= 4$, every vertex of $T$ is adjacent to at most two vertices of $V(C_4)$, and thus $T=T_0\cup T_1\cup T_2$, where $T_2=T_{\{x_1, x_3\}}\cup T_{\{x_2, x_4\}}$. Again by $g(G)=4$, we have $E(T_{\{x_i, x_{i+2}\}}, T_{x_i})=\emptyset$, $G[T_{x_i}]$ and $G[T_{\{x_i, x_{i+2}\}}]$ are empty graphs. Hence $|T_{x_i}|\leq 3$ because $U_{1,1,1,1}^4\not\subset G$.

(ii)\ By (i), we have $G[T_{\{x_1, x_3\}}]$ and $G[T_{\{x_2, x_4\}}]$ are empty. Suppose to the contrary that there exist vertices $u$, $v$ such that $u\in T_{\{x_1, x_3\}}$, $v\in T_{\{x_2, x_4\}}$ and $u\sim v$. Then $G[\{x_1,x_2,x_3,x_4,u,v\}]\cong K_{3,3}$, a contradiction.

(iii)\ It follows by $Z_3 \not \subset G$ and $g(G)=4$.

(iv)\ It follows by $U_2^4\not\subset G$.

(v)\ Suppose to the contrary that there exists a vertex $u\in T_0$ with $d_{T_2}(u)=0$. By (iv), $u$ must be adjacent to a vertex of $T_0$. Then there exists an edge in $G[T_0]$ with one endpoint is adjacent to a vertex of $T_2$. This implies $Z_1\subset G$, a contradiction. Thus, $d_{T_2}(u)\geq 1$ for every $u\in T_0$.

In what follows, we assume that $d_{T_{\{x_i, x_{i+2}\}}}(v)\geq 1$ for $v\in T_0$. Suppose that $d_{T_{\{x_{i+1}, x_{i+3}\}}}(v)\geq 1$, then $X_5$ (see Figure~\ref{Fig6}) is a subgraph of $G$. It implies that $G$ contains $K_{3,3}$ as a minor, a contradiction.
Suppose now that $d_{T_{\{x_i, x_{i+2}\}}}(v)\geq 3$, then $X_6\subset G$. Hence, $G$ contains $K_{3,3}$ as a minor, which is also a contradiction.

If $d_{T_{\{x_i, x_{i+2}\}}}(v)= 2$, we are to show that $|T_2|=|T_{\{x_i, x_{i+2}\}}|=2$. Note that $T_2=T_{\{x_i, x_{i+2}\}}\cup T_{\{x_{i+1}, x_{i+3}\}}$. If $T_{\{x_{i+1}, x_{i+3}\}}\neq \emptyset$, then $U^4_2\subset G$ by (ii), a contradiction. Hence, $T_{\{x_{i+1}, x_{i+3}\}}=\emptyset$ and so, $|T_2|=|T_{\{x_i, x_{i+2}\}}|$. Now we are to show $|T_{\{x_i, x_{i+2}\}}|=2$. 
Suppose to the contrary that $|T_{\{x_i, x_{i+2}\}}|\geq 3$. Then $Z_4\subset G$ by (ii), a contradiction. 

(vi)\ Suppose that there exist vertices $u$, $v$ of $T_0$ such that $u\sim v$.
By (v), we may assume that $u\sim u_1$ and $v\sim v_1$ with $u_1\in T_{\{x_s, x_{s+2}\}}$ and $v_1\in T_{\{x_t, x_{t+2}\}}$. By (iii), one has $u_1\neq v_1$, $u_1\not\sim v$ and $v_1\not \sim u$.
Then $G[\{x_{s-1},x_s,x_{s+1},u_1,u,v\}]\cong Z_1$, a contradiction.

(vii)\ Combining (i)  with $Z_1, Z_2\not \subset G$ proves (vii).

Our proof is complete.
\end{proof}

\begin{figure}[!ht]
\begin{center}
\includegraphics[scale=1]{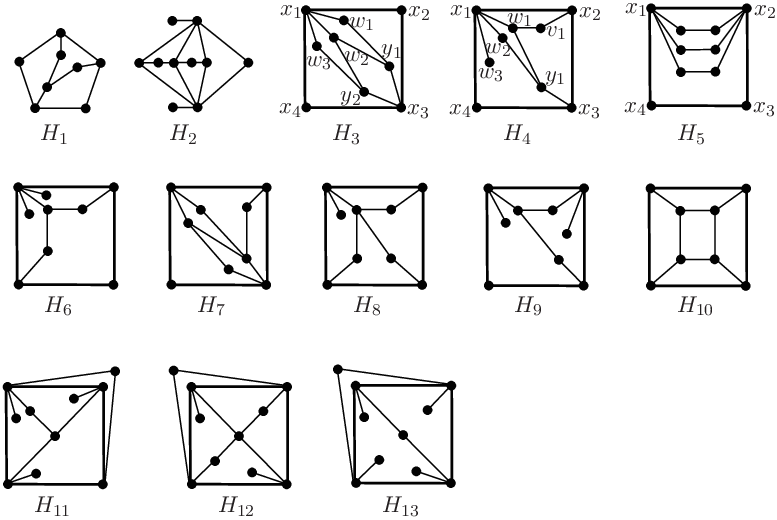}
\caption{Graphs $H_1,\ldots, H_{13}$.}\label{Fig3}
\end{center}
\end{figure}

Let $H_1,\ldots,H_{13}$ be the graphs illustrated in Figure~\ref{Fig3}. 

\begin{lemma}\label{lc4}
Let $G\in \mathcal{G}$ with girth $g(G)=4$, and let $C_4=x_1x_2x_3x_4x_1$ be its base graph. 
If $T_2=\emptyset$ and $\max\{|T_{x_i}|: 1\leq i \leq 4\}=3$, then $G$ is an induced subgraph of some graphs in $ \{H_3,H_4,H_5,H_6\}$. 
\end{lemma}

\begin{figure}[!ht]
\begin{center}
\includegraphics[scale=1.0]{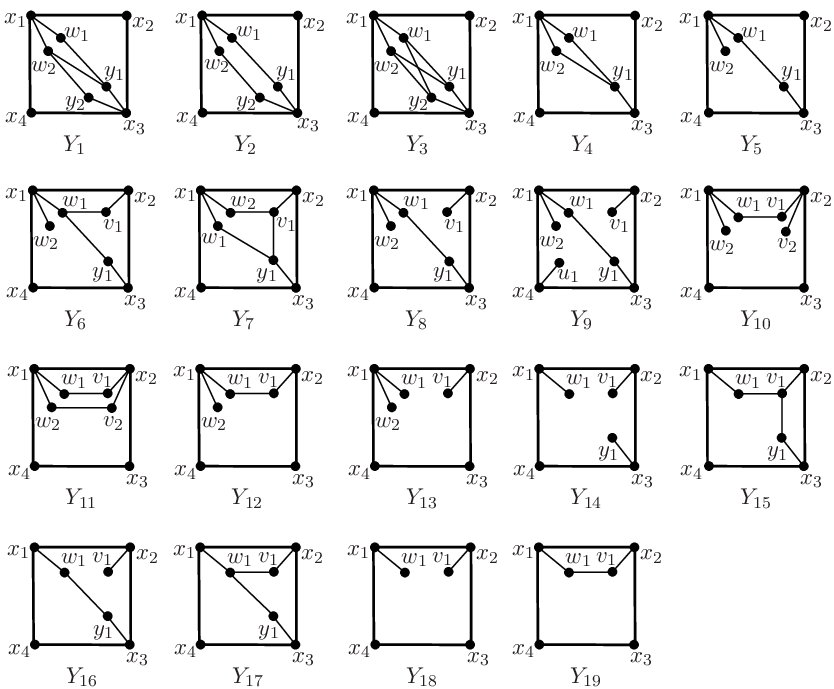}
\caption{Graphs $Y_1, Y_2, \ldots, Y_{19}$.}\label{Fig8}
\end{center}
\end{figure}

\begin{proof}
Combining $T_2=\emptyset$ with Lemma~\ref{c2}(vi) gives $T_0=\emptyset$. 
We assume, without loss of generality, that $|T_{x_1}|=\max\{|T_{x_i}|: 1\leq i \leq 4\}$ and $|T_{x_2}|\geq |T_{x_4}|$.  
Let $T_{x_1}=\{w_1, w_2, w_3\}$, $T_{x_2}=\{v_1, \ldots, v_{|T_{x_2}|}\}$, $T_{x_3}=\{y_1, \ldots, y_{|T_{x_3}|}\}$ and $T_{x_4}=\{u_1, \ldots, u_{|T_{x_4}|}\}$ (if $T_{x_i}\neq \emptyset$ for $2\leq i\leq 4$). 

We assert that $|T_{x_3}|\leq 2$. Suppose to the contrary that $|T_{x_3}|=3$. By Lemma~\ref{c2}(vii), we have $d_{T_{x_1}}(y_i)= 2$ and $d_{T_{x_3}}(w_i)= 2$, for $i\in\{1,2,3\}$. This implies that $X_7\subset G$. However, $X_7$ contains $K_{3,3}$ as a minor, which is impossible. 
Hence $|T_{x_3}|\le 2$, and so we proceed by considering the following three cases.

{\bf Case 1.} $|T_{x_3}|=2$. 

From Lemma~\ref{c2}(vii) we have
$H_3\subset G$.
We assert that $T_{x_2}=\emptyset$, then $T_{x_4}=\emptyset$ as $|T_{x_4}|\leq |T_{x_2}|$.
Suppose to the contrary that there exists a vertex, say $v_1$, in $T_{x_2}$.
From Lemma~\ref{c2}(vii), we get $d_{T_{x_1}}(v_1)=1$.
Firstly, we assume that $v_1\sim w_1$ (as the case of $v_1\sim w_3$ can be proved by symmetry, see Figure~\ref{Fig3}). We claim that $v_1\sim y_2$. Otherwise, $G[\{x_1,w_1,v_1,x_2,x_3,y_2\}]\cong U^4_2$, a contradiction.
Then $G$ contains $K_{3,3}$ as a minor, which contradicts the fact that $G$ is planar. 
Secondly, we assume that $v_1\sim w_2$.
Then $G$ contains $K_{3,3}$ as a minor, which is also a contradiction.
So we do indeed have $T_{x_2}=T_{x_4}=\emptyset$.
Hence
$G\cong H_3$.

{\bf Case 2.} $|T_{x_3}|=1$. 

Note $T_{x_3}=\{y_1\}$. By Lemma~\ref{c2}(vii), one has $d_{T_{x_1}}(y_1)=2$. Then we may assume that $N_{T_{x_1}}(y_1)=\{w_1,w_2\}$.

We first claim that $T_{x_4}=\emptyset$.
Indeed, if there exists a vertex $u_1\in T_{x_4}$, since $|T_{x_2}|\geq |T_{x_4}|$, there exists a vertex, say  $v_1$, in $T_{x_2}$.
By Lemma~\ref{c2}(vii), one has $d_{T_{x_1}}(v_1)=1$ and $d_{T_{x_1}}(u_1)=1$.
Firstly, we assume that $N_{T_{x_1}}(v_1)=N_{T_{x_1}}(u_1)=\{ w_i\}$ for some $i\in\{1,2\}$.
Then $G$ contains $K_{3,3}$ as a minor, a contradiction.
Secondly, we assume that $N_{T_{x_1}}(v_1)=N_{T_{x_1}}(u_1)=\{ w_3\}$.
Since $U^4_2\not \subset G$,
$v_1\sim y_1$ and $u_1\sim y_1$,  
which also contradicts the fact that $G$ is $K_{3,3}$-minor free.
Thirdly, we assume that $N_{T_{x_1}}(v_1)=\{ w_1\}$ and $N_{T_{x_1}}(u_1)=\{ w_2\}$.
Then we get $v_1\sim u_1$, as $U^4_2\not\subset G$. This implies that $G$ contains $K_{3,3}$ as a minor, a contradiction.
At last we assume that $N_{T_{x_1}}(v_1)=\{ w_1\}$ and $N_{T_{x_1}}(u_1)=\{ w_3\}$. By $U^4_2\not\subset G$, we get $v_1\sim u_1$ and $y_1\sim u_1$, which also contradicts the fact that $G$ is $K_{3,3}$-minor free.
This implies that $T_{x_4}=\emptyset$.

Now we assert that $|T_{x_2}|\leq 1$. Suppose to the contrary that $|T_{x_2}|\geq 2$, then there exist two vertices $v_1, v_2 \in T_{x_2}$.
From Lemma~\ref{c2}(vii), we get $d_{T_{x_1}}(v_1)=1$, $d_{T_{x_1}}(v_2)=1$ and $N_{T_{x_1}}(v_1)\neq N_{T_{x_1}}(v_2)$.
Firstly, we assume that $N_{T_{x_1}}(v_1)=\{w_1\}$ and $N_{T_{x_1}}(v_2)=\{w_2\}$, 
then $G$ contains $K_{3,3}$ as a minor, a contradiction.
Secondly, we assume that $N_{T_{x_1}}(v_1)=\{w_1\}$ and $N_{T_{x_1}}(v_2)=\{w_3\}$, 
then $v_2\sim y_1$ as $U^4_2\not\subset G$, and so $G$ contains $K_{3,3}$ as a minor, a contradiction.
This implies that $|T_{x_2}|\in\{0, 1\}$.

If $|T_{x_2}|=0$, then $G\subset H_4$, as desired. Then we assume that there exists exactly one vertex, say $v_1$, in $T_{x_2}$. By Lemma~\ref{c2}(vii), one has $d_{T_{x_1}}(v_1)=1$.
If $N_{T_{x_1}}(v_1)=\{w_1\}$ (as the case of $N_{T_{x_1}}(v_1)=\{w_2\}$ can be proved by symmetry), then
$G\cong H_{4}$, as desired.
If $N_{T_{x_1}}(v_1)=\{w_3\}$, then $v_1\sim y_1$ as $U^4_2\not\subset G$. Thus
$G\cong X_8$, one sees $X_8\cong H_3$.
Consequently, either $G\subset H_{4}$ or $G\subset H_{3}$.

{\bf Case 3.} $|T_{x_3}|=0$. 

Note that $|T_{x_2}|\geq |T_{x_4}|$. We claim that either $|T_{x_2}|\leq 1$ or $T_{x_4}= \emptyset$.
Suppose to the contrary that $|T_{x_2}|\geq 2$ and $|T_{x_4}|\geq 1$.
Using a similar discussion as before, 
we deduce that $G$ contains $K_{3,3}$ as a minor, a contradiction.

If $|T_{x_2}|\leq 1$, from Lemma~\ref{c2}(vii), we have $G\subset H_6$ or $G\subset X_9\cong H_4$.

If $T_{x_4}= \emptyset$, then $G\subset H_{5}$ by Lemma~\ref{c2}(vii).

This completes the proof of Lemma~\ref{lc4}.
\end{proof}
Let $Y_1, Y_2, \ldots, Y_{19}$ be the graphs illustrated in Figure~\ref{Fig8}.
\begin{lemma}\label{lc5}
Let $G\in \mathcal{G}$ with girth $g(G)=4$, and let $C_4=x_1x_2x_3x_4x_1$ be its base graph. 
If $T_2=\emptyset$ and $\max\{|T_{x_i}|: 1\leq i \leq 4\}=2$, then $G$ is an induced subgraph of some graphs in $ \mathcal{F}_1\cup \{H_3,H_4,H_5,H_7,H_8,H_9\}$.
\end{lemma}

\begin{proof}
We assume, without loss of generality, that $|T_{x_1}|=\max\{|T_{x_i}|: 1\leq i \leq 4\}$ and $|T_{x_2}|\geq |T_{x_4}|$.  
Let $T_{x_1}=\{w_1, w_2\}$, $T_{x_2}=\{v_1, \ldots, v_{|T_{x_2}|}\}$, $T_{x_3}=\{y_1, \ldots, y_{|T_{x_3}|}\}$ and $T_{x_4}=\{u_1, \ldots, u_{|T_{x_4}|}\}$ (if $T_{x_i}\neq \emptyset$ for $2\leq i\leq 4$). 
Note that $0\le |T_{x_3}|\le 2$. Hence, we proceed by considering the following three cases.

{\bf Case 1.} $|T_{x_3}|=2$. 
Then $T_{x_3}=\{y_1, y_2\}$.
By Lemma~\ref{c2}(vii), we have $1\leq d_{T_{x_3}}(w_i)\leq 2$ and $1\leq d_{T_{x_1}}(y_i)\leq 2$ for $i\in\{1,2\}$.
Then either $Y_1\subset G$ or $Y_2\subset G$ (because $Y_3$ contains $K_{3,3}$ as a minor).

{\bf Subcase 1.1.} $Y_1\subset G$.

$\bullet$ $T_{x_2}\cup T_{x_4}= \emptyset$. Then $G\cong Y_1$, and so
$G\subset H_{7}$, as desired.

$\bullet$ $T_{x_2}\cup T_{x_4}\neq \emptyset$.
Note that $|T_{x_2}|\geq |T_{x_4}|$. Then we assume that there exists a vertex $v$ in $T_{x_2}$. By Lemma~\ref{c2}(vii), one has $0\leq d_{T_{x_i}}(v)\leq 1$ for $i\in\{1,3\}$. We are to show that
$v$ is non-adjacent to any vertex of $\{w_1,y_2\}$.
Suppose that $v$ is adjacent to one vertex of $\{w_1,y_2\}$, without loss of generality, we assume that $v\sim w_1$. Then $v\sim y_2$ as $U^4_2\not\subset G$. Hence $G$ contains $K_{3,3}$ as a minor, a contradiction. This implies that $v\not\sim w_1$ and $v\not\sim y_2$. Combining this with $U^4_2\not\subset G$ and $g(G)=4$ gives us that
$v$ is adjacent to exactly one vertex of $\{w_2,y_1\}$. Then
$H_{7}\subset G$.
By a similar discussion as before, we get $T_{x_2}\cup T_{x_4} =\{v\}$. Therefore,
$G\cong H_{7}$, as desired.

{\bf Subcase 1.2.} $Y_2\subset G$.

$\bullet$  $T_{x_2}\cup T_{x_4}= \emptyset$. Then $G\cong Y_2$, and so $G\subset \mathcal{F}_1$, as desired.

$\bullet$  $T_{x_2}\cup T_{x_4}\neq \emptyset$.
Then there exists a vertex $v$ in $T_{x_2}\cup T_{x_4}$. We assert that $v$ is non-adjacent to any vertex of $\{w_1,w_2, y_1,y_2\}$. Without loss of generality, we assume that $v\in T_{x_2}$ (as the case of $v\in T_{x_4}$ can be prove by symmetry). Suppose that $v\sim w_1$, then $v\sim y_2$ as $U^4_2\not\subset G$. Hence, $G$ contains $K_{3,3}$ as a minor, a contradiction.
This implies that any vertex of $T_{x_2}\cup T_{x_4}$ is non-adjacent to any vertex of $\{w_1,w_2, y_1,y_2\}$. Now combine with $Z_1\not\subset G$, one has $|T_{x_2}|\leq 1$ and $|T_{x_4}|\leq 1$. Hence,
$G\subset W_{7}\subset \mathcal{F}_1$ (see Figure~\ref{Fig5}), as desired.

{\bf Case 2.} $|T_{x_3}|=1$.
By Lemma~\ref{c2}(vii), we have $1\leq d_{T_{x_1}}(y_1)\leq 2$. Then either $Y_4\subset G$ or $Y_5\subset G$.

{\bf Subcase 2.1.} $Y_4\subset G$.

$\bullet$ $T_{x_2}= \emptyset$. Then $G\cong Y_4$, and so
$G\subset
X_{10}\cong H_7$, as desired.

$\bullet$ $T_{x_2}\neq \emptyset$. Thus $1\leq|T_{x_2}|\leq 2$. Since $U^4_2\not\subset G$ and $g(G)=4$, we have
\begin{align}\label{et7}
\text{
$v$ is adjacent to exactly one vertex of $\{w_1,w_2,y_1\}$, for each $v\in T_{x_2}$.
}
\end{align}
Since $|T_{x_2}|\geq 1$, we may assume that there exists a vertex, say $v_1$, in $T_{x_2}$.

We first consider $v_1\sim w_1$ (as the case of $v_1\sim w_2$ can be prove by symmetry), and so $v_1$ is non-adjacent to any vertex of $\{w_2,y_1\}$ by \eqref{et7}.
Suppose that there exists a vertex $u_1$ in $T_{x_4}$.
Since $U^4_2\not\subset G$ and $g(G)=4$, $u_1$ is adjacent to one vertex of $\{w_1,v_1\}$, which contradicts our assumption that $G$ is planar. 
This implies that $T_{x_4}=\emptyset$. If $|T_{x_2}|= 1$, then
$G\subset
X_{10}\cong H_7$, as desired.
If $|T_{x_2}|= 2$, then there exists a vertex, say $v_2$, in $T_{x_2}\setminus \{v_1\}$. Combine \eqref{et7} with $U^4_2\not\subset G$ and $G$ is $K_{3,3}$-minor free, we have $v_2\sim y_1$, and so
$G\cong
X_{10}\cong H_7$, as desired.

We now consider the case of $v_1\sim y_1$. If $T_{x_4}=\emptyset$, combine \eqref{et7} with $|T_{x_2}|\leq 2$ and $Z_2\not\subset G$, we have
$G\subset
X_{10}\cong H_7$, as desired.
Thus we assume that $T_{x_4}\neq\emptyset$, and let $u_1\in T_{x_4}$.
We will see that $u_1\sim y_1$ and $u_1$ is non-adjacent to any vertex of $\{v_1,w_1,w_2\}$.
In fact, by $U^4_2\not\subset G$, one sees that $u_1$ is adjacent to at least one vertex of $\{v_1,y_1\}$
(Otherwise, $G[\{v_1,y_1,x_2,x_3,x_4,u_1\}]\cong U^4_2$, a contradiction).
If $u_1\sim v_1$, then $G$ contains $K_{3,3}$ as a minor, a contradiction. This means that $u_1\sim y_1$ and $u_1\not\sim v_1$, and hence $u_1$ is non-adjacent to any vertex of $\{w_1,w_2\}$, as $g(G)=4$. It follows that
$X_
{11}\subset G$.
Suppose that there exists a vertex $v_2\in T_{x_2}\setminus \{v_1\}$. By a similar discussion as before, we get $v_2$ is adjacent to exactly one vertex of $\{w_1,w_2\}$. Without loss of generality, we assume that $v_2\sim w_1$. Now by $U^4_2\not\subset G$, we have $u_1\sim v_2$ (Otherwise, $G[\{x_1,x_2,w_1,v_2,x_4,u_1\}]\cong U^4_2$, a contradiction). Hence $G$ contains $K_{3,3}$ as a minor, a contradiction. This implies that $|T_{x_2}|=1$, and so $|T_{x_4}|= 1$. Consequently,
$G\cong X_
{11}\cong H_3$, as desired.

{\bf Subcase 2.2.} $Y_5\subset G$. 

If $T_{x_2}\neq \emptyset$, then $1\leq |T_{x_2}|\leq 2$, and we let $v_1\in T_{x_2}$. We claim that
\begin{align}\label{et8}
\text{
$Y_i\subset G$ for some $i\in\{6,7,8\}$.
}
\end{align}
If $v_1\sim w_1$, then $v_1\not\sim y_1$, as $g(G)=4$. Hence $v_1\not\sim w_2$  (Otherwise, $G[\{x_1,w_1,w_2,v_1,y_1,x_3\}]\cong U^4_2$, a contradiction).
Thus $Y_6\subset G$.
If $v_1\sim w_2$, then $v_1\sim y_1$ (Otherwise, $G[\{x_1,x_2,v_1,w_2,x_3,y_1\}]\cong U^4_2$, a contradiction). Hence $v_1\not\sim w_1$, because $g(G)=4$.
It follows that $Y_7\subset G$.
If $v_1$ is non-adjacent to any vertex of $\{w_1,w_2\}$, by $U^4_2\not\subset G$, one has $v_1\not\sim y_1$. This implies that $Y_8\subset G$. Thus \eqref{et8} holds.

{\bf Subsubcase 2.2.1.} $|T_{x_2}|=1$ and $T_{x_4}=\emptyset$. Then \eqref{et8} means that
$G\cong Y_i$ for some $i\in\{6,7,8\}$. Therefore, either
$G\subset H_{9}$ or
$G\subset X_
{12}\cong H_8$, as desired.

{\bf Subsubcase 2.2.2.} $|T_{x_2}|=1$ and $T_{x_4}\neq\emptyset$. Then $|T_{x_4}|=1$ and let $u_1\in T_{x_4}$.
By \eqref{et8}, we consider the following three cases.

$\bullet$  $Y_6\subset G$. Then $u_1$ is adjacent to one vertex of $\{w_1,v_1\}$ (Otherwise, $G[\{x_1,x_2,w_1,v_1,x_4,u_1\}]\cong U^4_2$, a contradiction).
If $u_1\sim w_1$, then
$G\cong H_{8}$, as desired.
If $u_1\sim v_1$, then $G$ contains $K_{3,3}$ as a minor, a contradiction.

$\bullet$  $Y_7\subset G$. Then by $U^4_2\not\subset G$, we have $u_1$ is adjacent to one vertex of $\{w_2,v_1\}$, and so $G$ contains $K_{3,3}$ as a minor, a contradiction.

$\bullet$  $Y_8\subset G$. We assert that $u_1$ is non-adjacent to any vertex of $\{w_1,w_2\}$. If not, suppose first that $u_1\sim w_1$, then by $U^4_2\not\subset G$, one has $u_1\sim v_1$. It follows that $G$ contains $K_{3,3}$ as a minor, a contradiction. Now we suppose that $u_1\sim w_2$, then by $U^4_2\not\subset G$, one has $u_1\sim v_1$ and $u_1\sim y_1$, then $G$ contains $K_{3,3}$ as a minor, which is also a contradiction.
This implies that $u_1$ is non-adjacent to any vertex of $\{w_1,w_2\}$, and thus $u_1$ is non-adjacent to any vertex of $\{v_1,y_1\}$, by $Z_1\not\subset G$. Hence, $G\cong Y_9$, and so $G\subset \mathcal{F}_1$, as desired.

{\bf Subsubcase 2.2.3.} $|T_{x_2}|=2$.
Then there exists a vertex $v_2\in T_{x_2}\setminus \{v_1\}$. We have the following useful fact:
\begin{align}\label{et9}
\text{$G[\{x_1,x_2,x_3,x_4,w_1,w_2,y_1,v_2\}]\cong Y_i$ for some $i\in\{6,7,8\}$.}
\end{align}
It can be proved using a similar discussion as the proof of \eqref{et8}, whose procedure is omitted here.

In view of \eqref{et8}, we consider the following three subcases.

$\bullet$  $Y_6\subset G$. By \eqref{et9},
we first assume that $G[\{x_1,x_2,x_3,x_4,w_1,w_2,y_1,v_2\}]\cong Y_6$. Together with Lemma~\ref{c2}(i), one has $G[\{w_1,v_1,v_2,x_1,x_4,x_3\}]\cong Z_1$, a contradiction.
Next, we consider the case of $G[\{x_1,x_2,x_3,x_4,w_1,w_2,y_1,v_2\}]\cong Y_7$, then $G$ contains $K_{3,3}$ as a minor, a contradiction.
Finally, we consider the case of $G[\{x_1,x_2,x_3,x_4,w_1,w_2,y_1,v_2\}]\cong Y_8$, then
$H_{9}\subset G$.
We assert that $T_{x_4}=\emptyset$.
If not, suppose that there exists a vertex $u_1\in T_{x_4}$. $U^4_2\not\subset G$ means that $u_1$ is adjacent to one vertex of $\{w_1,v_1\}$. If $u_1\sim w_1$, then $u_1\sim v_2$ (Otherwise, $G[\{x_1,w_1,u_1,x_4,x_2,v_2\}]\cong U^4_2$, a contradiction). Hence, $G$ contains $K_{3,3}$ as a minor, a contradiction.
It follows that $u_1\sim v_1$, and so $G$ contains $K_{3,3}$ as a minor, which is also a contradiction.
This implies that $T_{x_4}=\emptyset$. Hence,
$G\cong H_{9}$, as desired.

$\bullet$  $Y_7\subset G$. By a  similar discussion we may obtain that $G[\{x_1,x_2,x_3,x_4,w_1,w_2,y_1,v_2\}]\not\cong Y_i$ for $i\in\{6,7\}$.
Then $G[\{x_1,x_2,x_3,x_4,w_1,w_2,y_1,v_2\}]\cong Y_8$ by \eqref{et9}, and thus
$X_
{12}\subset G$.
We assert that $T_{x_4}=\emptyset$.
If not, suppose that there exists a vertex $u_1\in T_{x_4}$. It follows that $u_1$ is adjacent to one vertex of $\{w_2,v_1\}$ (Otherwise, $G[\{x_1,x_2,v_1,w_2,x_4,u_1\}]\cong U^4_2$, a contradiction).
However, whether $u_1\sim w_2$ or $u_1\sim v_1$, $G$ contains $K_{3,3}$ as a minor, a contradiction. This implies that $T_{x_4}=\emptyset$. Then
$G\cong X_
{12}\cong H_8$, as desired.

$\bullet$  $Y_8\subset G$.
If $G[\{x_1,x_2,x_3,x_4,w_1,w_2,y_1,v_2\}]\cong Y_8$, then $Z_2\subset G$, a contradiction. By \eqref{et9}, one gets $G[\{x_1,x_2,x_3,x_4,w_1,w_2,y_1,v_2\}]\cong Y_i$ for some $i\in\{6,7\}$. Using a similar discussion as before gives us either
$G\cong H_{9}$ or
$G\cong H_{8}$, as desired.

{\bf Subsubcase 2.2.4.} $T_{x_2}=\emptyset$. Then $G\cong Y_5$, and so
$G\subset H_{8}$, as desired.

{\bf Case 3.} $T_{x_3}=\emptyset$. In this case, $0\leq |T_{x_2}|\leq 2$. So we consider the following three subcases.

{\bf Subcase 3.1.} $|T_{x_2}|=2$.
By Lemma~\ref{c2}(vii) and $Z_2\not\subset G$, we get either $Y_{10}\subset G$ or $Y_{11}\subset G$.

$\bullet$  $Y_{10}\subset G$.
If $T_{x_4}=\emptyset$, then $G\cong Y_{10}$, and thus
$G\subset X_
{12}\cong H_8$, as desired.
Then we assume that $T_{x_4}\neq\emptyset$,
thus $1\leq |T_{x_4}|\leq 2$.
Since $U^4_2\not\subset G$,
\begin{align}\label{et10}
\text{
$u$ is adjacent to one vertex of $\{w_1,v_1\}$ for each $u\in T_{x_4}$.
}
\end{align}

If $|T_{x_4}|=1$, then $T_{x_4}=\{u_1\}$.
By \eqref{et10}, we first assume that $u_1\sim w_1$, then $u_1\sim v_2$ (Otherwise,  $G[\{x_1,x_4,u_1,w_1,x_2,v_2\}]\cong U^4_2$, a contradiction). This means that
$G\cong
X_{12}\cong H_8$, as desired.
We now consider the case $u_1\sim v_1$, then either $G\cong H_{9}$ (if $u_1$ is non-adjacent to any vertex of $\{w_2,v_2\}$) or $G\cong X_
{10}\cong H_7$ (otherwise), as desired.

Next we prove that $|T_{x_4}|\neq 2$.
Suppose that there exists a vertex $u_2\in T_{x_4}\setminus \{u_1\}$, the discussion as before (in the case $|T_{x_4}|=1$) means that
$$
\text{
$G[\{x_1,x_2,x_3,x_4,w_1,w_2,v_1,v_2,u_2\}]\cong B$ for some $B\in\{H_{9},X_
{12},X_{10}\}$.
}
$$
Since $G$ is $\{K_5,K_{3,3}\}$-minor free, the only possibility is  $G[\{x_1,x_2,x_3,x_4,w_1,w_2,v_1,v_2,u_j\}]
\cong H_{9}$ for $j\in\{1,2\}$. It follows that $G[\{x_1,x_2,w_1,w_2,x_4,u_1,u_2\}]\cong Z_2$, a contradiction.

$\bullet$  $Y_{11}\subset G$.
If $T_{x_4}=\emptyset$, then $G\cong Y_{11}$, and thus $G\subset
H_5$, as desired.

Now we show that $T_{x_4} = \emptyset$. Suppose that there exists a vertex $u_1\in T_{x_4}$. Since $U^4_2\not\subset G$, $u_1$ is adjacent to one vertex of $\{w_1,v_1\}$ and one vertex of $\{w_2,v_2\}$, which contradicts our assumption that $G$
is planar.

{\bf Subcase 3.2.} $|T_{x_2}|=1$.
Lemma~\ref{c2}(vii) gives either $Y_{12}\subset G$ or $Y_{13}\subset G$.

$\bullet$  $Y_{12}\subset G$.
If $T_{x_4}=\emptyset$, then $G\cong Y_{12}$, and thus $G\subset
H_{9}$, as desired.
Now we assume that $T_{x_4}\neq\emptyset$, then $T_{x_4}=\{u_1\}$ as $|T_{x_4}|\leq |T_{x_2}|$.
Since $U^4_2\not\subset G$ and $g(G)=4$,
\begin{align*}
\text{
$u_1$ is adjacent to exactly one vertex of $\{w_1,v_1\}$.
}
\end{align*}
If $u_1\sim w_1$, then $u_1\not\sim w_2$ (Otherwise, $G[\{w_1,w_2,u_1,x_4,x_3,x_2\}]\cong Z_1$, a contradiction). It follows that $G\subset
H_{8}$, as desired.

If $u_1\sim v_1$, then $G\subset
X_9\cong H_4$, as desired.

$\bullet$  $Y_{13}\subset G$.
If $T_{x_4}=\emptyset$, then $G\cong Y_{13}$, and thus $G\subset
X_
{11}\cong H_3$, as desired.
Now we assume that $T_{x_4}\neq\emptyset$, then $T_{x_4}=\{u_1\}$.
By Lemma~\ref{c2}(vii), $0\leq d_{T_{x_1}}(u_1)\leq 1$.

If $d_{T_{x_1}}(u_1)=0$, then $u_1\not\sim v_1$ (Otherwise, $G[\{x_1,w_1,w_2,x_2,v_1,u_1\}]\cong Z_1$, a contradiction). It follows that $G\subset
X_
{11}\cong H_3$, as desired.

If $d_{T_{x_1}}(u_1)= 1$, we have $u_1\sim v_1$, as $U^4_2\not\subset G$. It follows that $G\subset
H_{9}$, as desired.

{\bf Subcase 3.3.} $T_{x_2}=\emptyset$. Then $T_{x_4}=\emptyset$, and so $G\subset
X_
{11}\cong H_3$, as desired.

This completes the proof of Lemma~\ref{lc5}.
\end{proof}

\begin{lemma}\label{lcl1}
Let $G\in \mathcal{G}$ with girth $g(G)=4$, and let $C_4=x_1x_2x_3x_4x_1$ be its base graph. 
If $T_2=\emptyset$ and $\max\{|T_{x_i}|: 1\leq i \leq 4\}=1$, then $G$ is an induced subgraph of some graphs in
$\mathcal{F}_1\cup \{H_3, H_{8}, H_{10}\}$.
\end{lemma}

\begin{proof}
We assume, without loss of generality, that $|T_{x_1}|=\max\{|T_{x_i}|: 1\leq i \leq 4\}$ and $|T_{x_2}|\geq |T_{x_4}|$.  
Let $T_{x_1}=\{w_1\}$, $T_{x_2}=\{v_1\}$, $T_{x_3}=\{y_1\}$ and $T_{x_4}=\{u_1\}$ (if $T_{x_i}\neq \emptyset$ for $2\leq i\leq 4$).

{\bf Case 1.} $|T_{x_3}|=|T_{x_2}|=1$.
We assert that
\begin{align}\label{et12}
\text{
$Y_i\subset G$ for some $i\in\{14,15,16,17\}$.
}
\end{align}
In fact, if $w_1\not\sim y_1$, then either $Y_{14}\subset G$ or $Y_{15}\subset G$, as $U^4_2\not\subset G$.
If $w_1\sim y_1$, then either $Y_{16}\subset G$ or $Y_{17}\subset G$, as $g(G)=4$. Thus \eqref{et12} holds.

{\bf Subcase 1.1.} $Y_{14}\subset G$.
If $T_{x_4}=\emptyset$, then $G\cong Y_{14}$, and thus $G\subset
H_{8}$, as desired.
If $T_{x_4}\neq\emptyset$, then $T_{x_4}=\{u_1\}$.

If $u_1$ is adjacent to every vertex of $\{w_1,v_1,y_1\}$, then $G\subset
H_{8}$, as desired.

If $u_1$ is adjacent to exactly two vertex of $\{w_1,v_1,y_1\}$, one has $U^4_2\subset G$, a contradiction.

If $u_1$ is adjacent to exactly one vertex of $\{w_1,v_1,y_1\}$, then $u_1\sim v_1$, as $U^4_2\not\subset G$. This implies that $G\subset \mathcal{F}_1$, as desired.

If $u_1$ is non-adjacent to any vertex of $\{w_1,v_1,y_1\}$, then $G\subset
H_{8}$, as desired.

{\bf Subcase 1.2.} $Y_{15}\subset G$.
If $T_{x_4}=\emptyset$, then $G\cong Y_{15}$, and thus $G\subset
H_{10}$, as desired.
If $T_{x_4}\neq\emptyset$, then $T_{x_4}=\{u_1\}$.
Since $U^4_2\not\subset G$ and $g(G)=4$, $u_1$ is adjacent to exactly one vertex of $\{w_1,v_1\}$ (resp., $\{y_1,v_1\}$). This implies either
$G\subset
H_{8}$ (if $u_1\sim v_1$) or
$G\cong
H_{10}$ (otherwise), as desired.

{\bf Subcase 1.3.} $Y_{16}\subset G$.
If $T_{x_4}=\emptyset$, then $G\cong Y_{16}$, and thus $G\subset \mathcal{F}_{1}$, as desired.
If $T_{x_4}\neq\emptyset$, then $T_{x_4}=\{u_1\}$.
Since $g(G)=4$, $u_1$ is adjacent to at most one vertex of $\{w_1,y_1\}$.
If $u_1\sim w_1$ (as the case of $u_1\sim y_1$ can be prove by symmetry), then $u_1\sim v_1$, by $U^4_2\not\subset G$. Hence $G$ contains $K_{3,3}$ as a minor, which is a contradiction.
Thus
$u_1$ is non-adjacent to any vertex of $\{w_1,y_1\}$. By $Z_6\not\subset G$, we have $u_1\not\sim v_1$. This implies $G\subset \mathcal{F}_1$, as desired.

{\bf Subcase 1.4.} $Y_{17}\subset G$.
If $T_{x_4}=\emptyset$, then $G\cong Y_{17}$, and thus $G\subset
H_{8}$, as desired.
If $T_{x_4}\neq\emptyset$, then $T_{x_4}=\{u_1\}$.
Since $U^4_2\not\subset G$ and $g(G)=4$, $u_1$ is adjacent to exactly one vertex of $\{w_1,v_1\}$.
If $u_1\sim v_1$, then $G$ contains $K_{3,3}$ as a minor, a contradiction. It means that
$u_1\sim w_1$, then $u_1\not\sim y_1$, as $g(G)=4$. This implies that $G\subset
H_{8}$, as desired.

{\bf Case 2.} $|T_{x_3}|=1$ and $T_{x_2}=\emptyset$.
Then $T_{x_4}=\emptyset$, as $|T_{x_4}|\leq |T_{x_2}|$.
Thus $G\subset
H_{8}$, as desired.

{\bf Case 3.} $T_{x_3}=\emptyset$ and $|T_{x_2}|=1$. 
Then either $Y_{18}\subset G$ or $Y_{19}\subset G$.

{\bf Subcase 3.1.} $Y_{18}\subset G$.
If $T_{x_4}=\emptyset$,
then $G\cong Y_{18}$ and so $G\subset \mathcal{F}_1$, as desired.
If $T_{x_4}\neq\emptyset$, then $T_{x_4}=\{u_1\}$.
If $u_1$ is non-adjacent to any vertex of $\{w_1,v_1\}$, then $G\subset \mathcal{F}_1$, as desired.
If $u_1\sim w_1$, then $u_1\sim v_1$ (Otherwise, $G[\{w_1,u_1,x_4,x_1,x_2,v_1\}]\cong U^4_2$, a contradiction).
Hence, $G\subset
H_{8}$, as desired.
If $u_1\not\sim w_1$ and $u_1\sim v_1$,
then $G\subset \mathcal{F}_1$, as desired.

{\bf Subcase 3.2.} $Y_{19}\subset G$.
If $T_{x_4}=\emptyset$, then $G\cong Y_{19}$ and so $G\subset
X_8\cong H_3$, as desired.
If $T_{x_4}\neq\emptyset$, then $T_{x_4}=\{u_1\}$.
Since $U^4_2\not\subset G$ and $g(G)=4$, $u_1$ is adjacent to exactly one vertex of $\{w_1,v_1\}$, which gives $G\subset
X_8\cong H_3$, as desired.

{\bf Case 4.} $T_{x_3}=T_{x_2}=\emptyset$.
Then $T_{x_4}=\emptyset$, so $G\subset
H_{8}$, as desired.

This completes the proof of Lemma~\ref{lcl1}.
\end{proof}
\section{\normalsize Triangle-free planar graphs with $\lambda_2\le 1$}\label{s3} 
In this section, we characterize all the connected triangle-free planar graphs with $\lambda_2\le 1.$ 
Our main result reads as follows.
\begin{thm}\label{Th21}
For a connected triangle-free planar graph $G$, $\lambda_2(G)\leq 1$ holds if and only if $G$ is an induced subgraph of a graph belonging to
$
\mathcal{F}_1\cup \{H_1,H_2,\ldots, H_{13}\}.
$
\end{thm}

In order to give the proof of Theorem~\ref{Th21}, we need to characterize the structure of $G$. By Lemma~\ref{lemT}, we consider graphs $G$ in $\mathcal{G}$ with at least one cycle.
From Lemma~\ref{lem31}, $\lambda_2(P_6)>1$, so we only need to consider graphs $G$ in $\mathcal{G}$ with $4\leq g(G)\leq 6$.

\begin{lemma}\label{lemg6}
Let $G\in \mathcal{G}$ with girth $g(G)=6$. Then $G\cong C_6$. Hence $G\subset \mathcal{F}_1$.
\end{lemma}
\begin{proof}
We use $C_6$ as the base graph. Since $g(G)= 6$, every vertex of $T=V(G)\setminus V(C_6)$ is adjacent to at most one vertex of $V(C_6)$, and thus $T=T_0\cup T_1$. Suppose there exists a vertex $v$ in $T_1$, then $P_6\subset G$, a contradiction. Hence, $T_1=\emptyset$, and so $T_0=\emptyset$. This completes the proof.
\end{proof}

\begin{lemma}\label{lemg5}
Let $G\in \mathcal{G}$ with girth $g(G)=5$. Then $G$ is an induced subgraph of $H_1$ or some member in $\mathcal{F}_1$.
\end{lemma}
\begin{proof}
We use $C_5=x_1 x_2 x_3 x_4 x_5 x_1$ as the base graph. Since $g(G)= 5$, every vertex of $T$ is adjacent to at most one vertex of $V(C_5)$, and thus $T=T_0\cup T_1$.
Suppose $T_0\neq \emptyset$, then there exists a vertex of $T_0$ which is adjacent to at least one vertex of $T_1$. This implies that $P_6\subset G$, a contradiction. Hence, $T_0= \emptyset$. Consequently, $T=T_1=\bigcup_{i=1}^5 T_{x_i}$. By $Z_1\not\subset G$ and $g(G)=5$, we immediately get the following claim.
\begin{claim}\label{c1} Let $i$, $j$ be distinct integers of $\{1,2, \ldots, 5\}$. Then the following statements hold:
\begin{wst}
\item[{\rm (i)}] $|T_{x_i}| \leq 1$;
\item[{\rm (ii)}] If $u\in T_{x_i}$ and $v\in T_{x_j}$, then $u\not \sim v$ for $x_i \sim x_j$ and $u \sim v$ for $x_i \not \sim x_j$.
\end{wst}
\end{claim}


Now we come back to show our result. By Claim~\ref{c1}(i), we have $0\leq |T_1|\leq 5$.

If $0\leq |T_1|\leq 3$, using Claim~\ref{c1}(ii), we deduce either $G\subset  X_4 \subset \mathcal{F}_1$ 
or $G\subset H_1$.

If $|T_1|=4$, by Claim~\ref{c1}(ii), we get $G\cong X_3$. One sees that $X_3$ contains $K_{3,3}$ as a minor, a contradiction.

If $|T_1|=5$, then Claim~\ref{c1}(ii) gives $G\cong X_2$, i.e., $G$ is the Petersen graph. One sees that $X_2$ contains $K_5$ as a minor, a contradiction.

This completes the proof of Lemma~\ref{lemg5}.
\end{proof}

\begin{lemma}\label{lemg4}
Let $G\in \mathcal{G}$ with girth $g(G)=4$. Then $G$ is an induced subgraph of some graph of $\mathcal{F}_1\cup \{H_2,H_3,\ldots, H_{13}\}$.
\end{lemma}

\begin{proof}
We use $C_4=x_1 x_2 x_3 x_4 x_1$ as the base graph (in what follows, the subscripts are taken modulo 4). From Lemma~\ref{c2}(i), $T=T_0\cup T_1\cup T_2$, where $T_2=T_{\{x_1, x_3\}}\cup T_{\{x_2, x_4\}}$. 
We prove our result by distinguishing the following four cases depending on the structure of $T$.

{\bf Case 1.} $T_1=\emptyset$ and $T_2\neq\emptyset$. 

Without loss of generality, we may assume that $|T_{\{x_1, x_3\}}|\geq |T_{\{x_2, x_4\}}|$. We assert that $|T_{\{x_2, x_4\}}|\leq 1$. If not, suppose $|T_{\{x_2, x_4\}}|\geq 2$. Together with Lemma~\ref{c2}(ii), we have $Z_3\subset G$, a contradiction.
So we do indeed have $|T_{\{x_2, x_4\}}|\leq 1$.

{\bf Subcase 1.1.} $|T_{\{x_2, x_4\}}|= 0$. 
Then $T_2=T_{\{x_1, x_3\}}$.
If $T_0=\emptyset$, then $G\subset \mathcal{F}_1$. Hence, we assume that $T_0\neq \emptyset$.

If there exists a vertex $v$ of $T_0$ such that $d_{T_2}(v)=2$, by Lemma~\ref{c2}(v), we have $|T_2|= 2$. Combine with Lemma~\ref{c2}(iii), we get
$G\subset H_2$, as desired.
Now by Lemma~\ref{c2}(v), we only need to consider that $d_{T_2}(v)= 1$ for every vertex $v$ of $T_0$. In this subcase, combining with Lemma~\ref{c2}(iii) and (vi), we have $G\subset \mathcal{F}_1$, as desired.

{\bf Subcase 1.2.} $|T_{\{x_2, x_4\}}|= 1$. We may assume that $T_{\{x_2, x_4\}}=\{w\}$.
By Lemma~\ref{c2}(ii) and $Z_4\not\subset G$, we obtain $|T_{\{x_1, x_3\}}|\leq 2$.
We are to show that $T_0=\emptyset$.
Suppose that there exists a vertex $v$ of $T_0$. If $v\sim w$, since $|T_{\{x_1, x_3\}}|\geq |T_{\{x_2, x_4\}}|$, there exists a vertex, say $u$, in $T_{\{x_1, x_3\}}$.
By Lemma~\ref{c2}(ii) and (v), we have $u\not\sim w$ and $v\not\sim u$, and so $G[\{x_1,x_2,x_3,u,w,v\}]\cong U^4_2$, a contradiction.
This implies that $v\not\sim w$.
It is similarly to show that $v$ is also non-adjacent to any vertex of $T_{\{x_1, x_3\}}$. This contradicts Lemma~\ref{c2}(vi). Hence $T_0=\emptyset$, and so $G\subset H_2$, as desired.

{\bf Case 2.} $T_1\neq \emptyset$ and $T_2=\emptyset$.

Lemma~\ref{c2}(vi) gives us $T_0=\emptyset$. 
By $T_1\neq\emptyset$ and Lemma~\ref{c2}(i), we have $1\le \max\{|T_{x_i}|: 1\leq i \leq 4\}\le 3$.
If $\min\{|T_{x_i}|:1\leq i\leq 4\}=3$, from Lemma~\ref{c2}(vii) and $g(G)=4$, we have $G\cong X_1$.
However, one sees that $X_2\subset X_1$, and so $X_1$ contains $K_5$ as a minor, a contradiction. Thus, $\min\{|T_{x_i}|:1\leq i\leq 4\}\le 2$, and then from Lemmas~\ref{lc4}-\ref{lcl1}, we get that 
$G$ is an induced subgraph of some graphs in $ \mathcal{F}_1\cup \{H_3,\ldots,H_{10}\}$.

%
%
%
%

{\bf Case 3.} $T_1\neq \emptyset$ and $T_2\neq \emptyset$.

Since $G$ is $K_{3,3}$-minor free, we immediately get the following fact.
\begin{fact}\label{f0}
For $i=1,2,3,4,$ every vertex of $T_{x_i}$ is adjacent to at most one vertex of $T_{\{x_{i+1},x_{i+3}\}}.$
\end{fact}

We need the following claim.
\begin{claim}\label{c6}
Assume $i$ is an integer in $\{1,2,3,4\}.$ 
\begin{wst}
\item[{\rm (i)}] If $e(T_{x_i}, T_{x_{i+1}})\geq 1$, then $|T_2|= 1$;
\item[{\rm (ii)}] If $e(T_{x_i}, T_{x_{i+2}})\geq 1$, then $T_{\{x_{i+1}, x_{i+3}\}}=\emptyset$, hence $T_2=T_{\{x_i, x_{i+2}\}}.$ 
\end{wst}
\end{claim}
\begin{proof}[\bf Proof of Claim~\ref{c6}]
(i)\ Since $T_2\neq\emptyset$, we suppose to the contrary that $|T_2|\geq 2$. Then there exist two vertices $w_1,w_2\in T_2$.
By $e(T_{x_i}, T_{x_{i+1}})\geq 1$, we may assume that there exist vertices $u,v$ such that $u\in T_{x_i}$, $v\in T_{x_{i+1}}$ and $u\sim v$.

Firstly, we assume that $w_1,w_2\in T_{\{x_i, x_{i+2}\}}$ (as the case of $w_1,w_2\in T_{\{x_{i+1}, x_{i+3}\}}$ can be proved by symmetry).
By $g(G)=4$ and $U^4_2\not\subset G$, we have $w_j\not\sim u$ and $w_j\sim v$ for $j\in\{1,2\}$, contradicting Fact~\ref{f0}. 

Next, we consider the case of $w_1\in T_{\{x_i, x_{i+2}\}}$ and $w_2\in T_{\{x_{i+1}, x_{i+3}\}}$. By a similar discussion as before, we get $w_1\sim v$ and $w_2\sim u$, which contradicts the fact that $G$ is $K_{3,3}$-minor free.

(ii)\ Suppose to the contrary that there exists a vertex $w\in T_{\{x_{i+1}, x_{i+3}\}}$.
By $e(T_{x_i}, T_{x_{i+2}})\geq 1$, we may assume that there exist vertices $u,v$ such that $u\in T_{x_i}$, $v\in T_{x_{i+2}}$ and $u\sim v$.
Then $w$ is adjacent to at least one vertex of $\{u,v\}$, as $U^4_2\not\subset G$.
Hence, $G$ contains $K_{3,3}$ as a minor, a contradiction.

This completes the proof of Claim~\ref{c6}.
\end{proof}

Lemma~\ref{c2}(i) implies that $G[T_{x_i}]$ is empty for each $i\in\{ 1,2,3,4\}$. Then we proceed by considering the subsequent three subcases.

{\bf Subcase 3.1.} $T_1$ is an independent set of $G$.

Here, we assume that $|T_{x_1}|=\max\{|T_{x_i}|: 1\leq i \leq 4\}$ and $|T_{x_2}|\geq |T_{x_4}|$. By $T_1\neq \emptyset$ and Lemma~\ref{c2}(i), $1\leq|T_{x_1}|\leq 3$.
Let $W_1,W_2,\ldots, W_{7}$ be the graphs illustrated in Figure~\ref{Fig5}.
In order to show our result, we need the subsequent claims.
\begin{claim}\label{c7}
If $|T_{x_1}|=3$, then
$G\subset W_{1} \subset
H_2$.
\end{claim}
\begin{proof}[\bf Proof of Claim~\ref{c7}]
If $|T_{x_1}|=3$, by Lemma~\ref{c2}(vii) and $E(G[T_1])=\emptyset$, we have $T_{x_i}=\emptyset$ for $i \in\{2,3, 4\}$.

We now show that $T_2=T_{\{x_2,x_4\}}$. Suppose to the contrary that there exists a vertex, say $v$, in $T_{\{x_1,x_3\}}$.
It implies that $Z_5\subset G$ by Lemma~\ref{c2}(i), a contradiction.
Hence $|T_{\{x_2,x_4\}}|\geq 1$, as $T_2\neq \emptyset$.

Then we show that $|T_{\{x_2,x_4\}}|=1$. Suppose to the contrary that there exist two vertices $u_1, u_2$ in $T_{\{x_2,x_4\}}$.
By $Z_3,Z_4\not\subset G$, we have $d_{T_{x_1}}(u_i)= 2$ for $i\in\{1, 2\}$, which contradicts Fact~\ref{f0}. 
Hence $|T_2|=|T_{\{x_2,x_4\}}|=1$,
and the vertex of $T_{\{x_2,x_4\}}$ is adjacent to exactly two vertices of $T_{x_1}$. Combine with Lemma~\ref{c2}(vi), one has $G\subset W_1$, as desired.
\end{proof}
\begin{figure}[!ht]
\begin{center}
\includegraphics[scale=0.9]{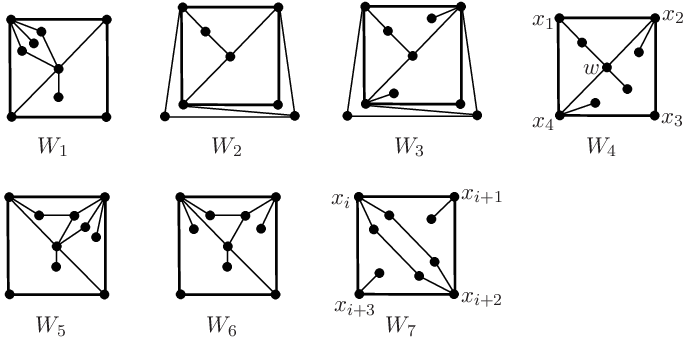}
\caption{Graphs $W_1, W_2, \ldots, W_{7}$.}\label{Fig5}
\end{center}
\end{figure}

If $|T_{x_1}|=2$, together with $E(G[T_1])=\emptyset$ and $Z_1, Z_2 \not\subset G$, one has
\begin{align}\label{et1}
\text{
$T_{x_3}=\emptyset$ and $\max\{|T_{x_2}|,|T_{x_4}|\}\leq 1$.
}
\end{align}
We proceed with the following two claims for the case $|T_{x_1}|=2$.

\begin{claim}\label{c8}
If  $|T_{x_1}|=2$ and $T_{\{x_2,x_4\}}\neq \emptyset$, then the following statements hold:
\begin{wst}
\item[{\rm (i)}]
If $v\in T_{\{x_2,x_4\}}$, then $d_{T_{x_1}}(v)\in\{1,2\};$
\item[{\rm (ii)}]
$1\leq|T_{\{x_2,x_4\}}|\leq 2;$
\item[{\rm (iii)}]
If there exists $v\in T_{\{x_2,x_4\}}$ with $d_{T_{x_1}}(v)=2$, then
$G\subset
H_2;$
\item[{\rm (iv)}]
If $|T_{\{x_2,x_4\}}|= 2$ and $d_{T_{x_1}}(v)=1$ for each $v\in T_{\{x_2,x_4\}}$, then
$G\subset
H_2;$
\item[{\rm (v)}]
If $T_{\{x_2,x_4\}}= \{v\}$ and $d_{T_{x_1}}(v)=1$, then
$G\subset
H_{11}$.
\end{wst}
\end{claim}
\begin{proof}[\bf Proof of Claim~\ref{c8}]
(i)\ Suppose that there exists a vertex $v\in T_{\{x_2,x_4\}}$ with $d_{T_{x_1}}(v)=0$, then $Z_3\subset G$, a contradiction.

(ii)\ Combining (i) with Fact~\ref{f0}, we get $|T_{\{x_2,x_4\}}|\leq 2$.

(iii)\ First, we assert that $T_{x_2}=T_{x_4}=\emptyset$.
Suppose that there exists a vertex in $T_{x_2}$ (resp., $T_{x_4}$). Combine $E(G[T_1])=\emptyset$ with Lemma~\ref{c2}(i), one has $Z_3\subset G$, a contradiction.
It follows that $T_{x_2}=T_{x_4}=\emptyset$.
Next we show that $T_{\{x_2,x_4\}}= \{v\}$.
Suppose that there exists a vertex $u\in T_{\{x_2,x_4\}}\setminus \{v\}$, then $d_{T_{x_1}}(u)\in\{1,2\}$ by (i). However, this contradicts Fact~\ref{f0} and hence we conclude that $T_{\{x_2,x_4\}}= \{v\}$.
We now claim that $T_{\{x_1,x_3\}}= \emptyset$. Indeed, if there exists a vertex $w\in T_{\{x_1,x_3\}}$, by Lemma~\ref{c2}(i)-(ii), we get $Z_1\subset G$, a contradiction. It implies that $T_{\{x_1,x_3\}}= \emptyset$.
Consequently, $T_{x_2}=T_{x_4}=T_{\{x_1,x_3\}}=\emptyset$ and $T_{\{x_2,x_4\}}=\{v\}$.
Together with Lemma~\ref{c2}(vi) and \eqref{et1}, we get $G\subset W_{1}\subset H_2$, as desired.

(iv)\ If $|T_{\{x_2,x_4\}}|= 2$, we assume that $T_{\{x_2,x_4\}}=\{v_1,v_2\}$. The assumption gives us $d_{T_{x_1}}(v_1)=1$ and $d_{T_{x_1}}(v_2)=1$. 
Then $N_{T_{x_1}}(v_1)\neq N_{T_{x_1}}(v_2)$ by Fact~\ref{f0}.
Hence, $W_2\subset G$.

We first prove $T_{\{x_1,x_3\}}=\emptyset$.
Suppose that there exists a vertex in $T_{\{x_1,x_3\}}$, then $G$ contains $K_{3,3}$ as a minor, a contradiction.

Next, we claim that $T_0=\emptyset$. Indeed, if there exists a vertex, say $w$, in $T_0$, then $d_G(w)=d_{T_{\{x_2,x_4\}}}(w)\in\{1,2\}$ by Lemma~\ref{c2}(v)-(vi).
If $d_{T_{\{x_2,x_4\}}}(w)=1$, then $Z_3\subset G$, which is a contradiction. If $d_{T_{\{x_2,x_4\}}}(w)=2$, then $G$ contains $K_5$ as a minor, a contradiction.
This implies that $T_{\{x_1,x_3\}}=T_0=\emptyset$.
Together with \eqref{et1} and Lemma~\ref{c2}(i), one has $G\subset
W_{3}\cong H_2$.

(v)\ We are to show that $|T_{\{x_1,x_3\}}|\leq 1$. Suppose to the contrary that $|T_{\{x_1,x_3\}}|\geq 2$. From Lemma~\ref{c2}(i)-(ii), one has $U^4_2\subset G$, a contradiction.
This implies that $|T_{\{x_1,x_3\}}|\leq 1$.

Next, we assert that $T_0=\emptyset$. Suppose that there exists a vertex, say $u$, in $T_0$.
By Lemma~\ref{c2}(v)-(vi), we have $d_{T_2}(u)=1$, and so $U^4_2\subset G$, a contradiction.
This implies that $T_0=\emptyset$.
Now, combine Lemma~\ref{c2}(ii) with \eqref{et1}, one has $G\subset
H_{11}$.

This completes the proof of Claim~\ref{c8}.
\end{proof}

\begin{claim}\label{cl8}
If  $|T_{x_1}|=2$ and $T_{\{x_2,x_4\}}= \emptyset$, then $G \subset \mathcal{F}_1$.
\end{claim}

\begin{proof}[\bf Proof of Claim~\ref{cl8}]
Since $T_2\neq \emptyset$, we have $T_{\{x_1,x_3\}}\neq \emptyset$.
We first assert that
\begin{align}\label{et2}
\text{
$d_G(v)=d_{T_{\{x_1,x_3\}}}(v)=1$ for each $v\in T_0$.
}
\end{align}
Combine $T_{\{x_2,x_4\}}= \emptyset$ with Lemma~\ref{c2}(v)-(vi), we have $d_G(v)=d_{T_{\{x_1,x_3\}}}(v)\in\{1,2\}$ for each $v\in T_0$.
Suppose that there exists a vertex $v_0\in T_0$ such that $d_{T_{\{x_1,x_3\}}}(v_0)=2$, then $Z_3\subset G$, a contradiction.
Thus \eqref{et2} holds.

Next, we assert that
\begin{align}\label{et3}
\text{
$E(T_{\{x_1,x_3\}},T_{x_i})=\emptyset$, for each $i\in\{2,4\}$.
}
\end{align}
Indeed, if there exists an edge in $E(T_{\{x_1,x_3\}},T_{x_i})$ for some $i\in\{2,4\}$, then $Z_3\subset G$, a contradiction. Hence \eqref{et3} holds.

Now, Combine \eqref{et1}-\eqref{et3} with Lemma~\ref{c2}(i), one has $G\subset \mathcal{F}_1$.
This completes the proof of Claim~\ref{cl8}.
\end{proof}

Further on, we need the following three claims for the case of $|T_{x_1}|=1$.
\begin{claim}\label{c9}
Assume  $|T_{x_1}|=1$, $T_{\{x_1,x_3\}}= \emptyset$ and $T_{\{x_2,x_4\}}\neq \emptyset.$
\begin{wst}
\item[{\rm (i)}] If there exist a vertex $w\in T_{\{x_2,x_4\}}$ and a number $s\in\{1,3\}$, such that $d_{T_{x_s}}(w)=1$, then
    $G\subset
    H_2;$
\item[{\rm (ii)}] If $N_{T_{x_1}}(w)=N_{T_{x_3}}(w)=\emptyset$ for each $w\in T_{\{x_2,x_4\}}$,
    then $G\subset \mathcal{F}_1.$
\end{wst}
\end{claim}

\begin{proof}[\bf Proof of Claim~\ref{c9}]
(i) Suppose that $|T_{x_t}|=1$ and $d_{T_{x_t}}(w)=0$ for $t\in\{1,3\}\setminus \{s\}$, then $U^4_2\subset G$, a contradiction. It implies that
\begin{align*}
\text{
either $T_{x_t}=\emptyset$ or $d_{T_{x_t}}(w)=1$,
}
\end{align*}
for $t\in\{1,3\}\setminus \{s\}$.

Firstly, we consider that $T_{x_t}=\emptyset$. Since $|T_{x_1}|=1$, we get $t=3$, and hence $s=1$. Combine Fact~\ref{f0} with Lemma~\ref{c2}(ii) and $Z_4\not\subset G$, one has $1\leq |T_{\{x_2,x_4\}}|\leq 2$.

$\bullet$ $|T_{\{x_2,x_4\}}|= 2$ and $T_0\neq \emptyset$. We assert that
\begin{align}\label{et4}
\text{
$d_{T_{\{x_2, x_4\}}}(u)=2$ for each $u\in T_0$, and so $|T_0|=1$.
}
\end{align}
In fact, if there exists a vertex $u\in T_0$, then $d_{T_{\{x_2, x_4\}}}(u)\in\{1,2\}$ by Lemma~\ref{c2}(v).
Suppose to the contrary that $d_{T_{\{x_2, x_4\}}}(u)=1$, then either $Z_3\subset G$ (if $u\sim w$) or $U^4_2\subset G$ (if $u$ is adjacent to the vertex of $T_{\{x_2,x_4\}} \setminus \{w\}$), both are impossible.
It implies that $d_{T_{\{x_2, x_4\}}}(u)=2$ for each $u\in T_0$.
Now by Lemma~\ref{c2}(iii), we have $|T_0|=1$. Hence, \eqref{et4} holds.
Combine \eqref{et4} with $|T_{\{x_2,x_4\}}|= 2$, $|T_{x_1}|= 1$, $T_{x_3}=\emptyset$, Fact~\ref{f0}, Lemma~\ref{c2}(i) and (vi), one has $G\subset
H_2$.

$\bullet$ $|T_{\{x_2,x_4\}}|= 2$ and $T_0=\emptyset$. Combining Fact~\ref{f0} with Lemma~\ref{c2}(i)-(ii) gives us $G\subset
H_2$.

$\bullet$ $|T_{\{x_2,x_4\}}|= 1$. Then $|T_0|\leq 1$ by Lemma~\ref{c2}(iii) and (vi).
Together with $T_{x_3}=\emptyset$, we obtain $G\subset W_4$, and so $G\subset
H_2$.

Secondly, we consider the case of $d_{T_{x_t}}(w)=1$. Then  $d_{T_{x_1}}(w)=d_{T_{x_3}}(w)=1$.
By Fact~\ref{f0}, Lemma~\ref{c2}(ii) and $Z_3\not\subset G$, we get
$1\leq|T_{\{x_2,x_4\}}|\leq 2$.

$\bullet$ $|T_{\{x_2,x_4\}}|= 2$.
We let $\{w^*\}=T_{\{x_2,x_4\}}\setminus \{w\}$.
Recall that $T_{\{x_1,x_3\}}=\emptyset$. Then $T_2=T_{\{x_2,x_4\}}=\{w,w^*\}$.
We claim that $T_0= \emptyset$.
Suppose to the contrary that there exists a vertex, say $v$, in $T_0$. By Lemma~\ref{c2}(v)-(vi), we have $d_G(v)=d_{T_{\{x_2,x_4\}}}(v)\in \{1,2\}$.

If $d_{T_{\{x_2,x_4\}}}(v)=1$ and $N_{T_{\{x_2,x_4\}}}(v)=\{w\}$, then $G[\{x_1,x_2,w,x_4,w^*,v\}\cup T_{x_3}]\cong Z_3$ by Fact~\ref{f0} and Lemma~\ref{c2}(ii), a contradiction.

If $d_{T_{\{x_2,x_4\}}}(v)=1$ and $N_{T_{\{x_2,x_4\}}}(v)= \{w^*\}$, then $G[T_{x_1}\cup\{x_1,x_2,w,w^*,v\}]\cong U^4_2$ by Fact~\ref{f0} and Lemma~\ref{c2}(ii), a contradiction.

If $d_{T_{\{x_2,x_4\}}}(v)=2$, then $G$ contains $K_{3,3}$ as a minor, a contradiction.

This implies that $T_0=\emptyset$. Together with $T_2=\{w,w^*\}$, Fact~\ref{f0} and Lemma~\ref{c2}(ii), one obtains that $G\subset
H_2$.

$\bullet$ $|T_{\{x_2,x_4\}}|= 1$. That is, $T_2=T_{\{x_2,x_4\}}= \{w\}$ as $T_{\{x_1,x_3\}}=\emptyset$.
Together with $d_{T_{x_1}}(w)=d_{T_{x_3}}(w)=1$, we obtain $G\subset
H_2$.

(ii) We assert that
\begin{align}\label{et5}
\text{
$d_G(v)=d_{T_{\{x_2,x_4\}}}(v)=1$ for each $v\in T_0$.
}
\end{align}
Indeed, by $T_{\{x_1,x_3\}}=\emptyset$ and Lemma~\ref{c2}(v)-(vi), we have $d_G(v)=d_{T_{\{x_2,x_4\}}}(v)\in\{1,2\}$ for each $v\in T_0$. Suppose to the contrary that there exists a vertex $v_0\in T_0$ such that $d_{T_{\{x_2,x_4\}}}(v_0)=2$.
Since $|T_{x_1}|=1$ and $N_{T_{x_1}}(w)=N_{T_{x_3}}(w)=\emptyset$ for each $w\in T_{\{x_2,x_4\}}$, we get $U^4_2\subset G$, a contradiction. Hence, \eqref{et5} holds.
Now, Combine \eqref{et5} with Lemma~\ref{c2}(i)-(iii), one has $G\subset \mathcal{F}_1$.

This completes the proof of Claim~\ref{c9}.
\end{proof}

\begin{claim}\label{c10}
Assume $|T_{x_1}|=1$, $T_{\{x_1,x_3\}}\neq\emptyset$ and $T_{\{x_2,x_4\}}= \emptyset.$ 
\begin{wst}
\item[{\rm (i)}] If there exist a vertex $w\in T_{\{x_1,x_3\}}$ and a number $s\in\{2,4\}$, such that $d_{T_{x_s}}(w)=1$, then 
    $G\subset
    H_2;$
\item[{\rm (ii)}] If $N_{T_{x_2}}(w)=N_{T_{x_4}}(w)=\emptyset$ for each $w\in T_{\{x_1,x_3\}}$,
    then $G\subset \mathcal{F}_1.$
\end{wst}
\end{claim}
We omit the proof of Claim~\ref{c10} due to the similarity with the proof of Claim~\ref{c9}.

\begin{claim}\label{c11}If $|T_{x_1}|=1$,
$\min\{|T_{\{x_1,x_3\}}|,|T_{\{x_2,x_4\}}|\}\geq 1$, then $G$ is an induced subgraph of some graph in $\{H_2,H_{12},H_{13}\}$.
\end{claim}

\begin{proof}[\bf Proof of Claim~\ref{c11}]
By Lemma~\ref{c2}(ii) and $Z_3\not\subset G$, we get
$\min\{|T_{\{x_1,x_3\}}|,|T_{\{x_2,x_4\}}|\}= 1$.
By Lemma~\ref{c2}(ii) and $Z_4\not\subset G$, one has $\max\{|T_{\{x_1,x_3\}}|,|T_{\{x_2,x_4\}}|\}\leq 2$.

If $|T_{\{x_1,x_3\}}|\leq |T_{\{x_2,x_4\}}|$. Then there exists exactly one vertex, say $z$, in $T_{\{x_1,x_3\}}$.

$\bullet$ $|T_{\{x_2,x_4\}}|= 2$. Since $Z_3\not\subset G$ and Fact~\ref{f0}, the vertex of $T_{x_1}$ (say $v$) is adjacent to exactly one vertex of $T_{\{x_2,x_4\}}$ (say $w$).
We first assert that $T_{x_3}=\emptyset$.
Suppose that there exists a vertex, say $u$, in $T_{x_3}$.
If $u\sim w$, then $G$ contains $K_5$ as a minor, a contradiction.
If $u\not\sim w$, then $G[\{x_1,v,w,x_4,x_3,u\}]\cong U^4_2$, a contradiction.
This implies that $T_{x_3}=\emptyset$.
Next we assert that $T_0=\emptyset$.
Suppose that there exists a vertex, say $y$, in $T_0$.
Combine $T_{\{x_1,x_3\}}=\{z\}$, $|T_{\{x_2,x_4\}}|= 2$ with Lemma~\ref{c2}(v)-(vi), one has either $N_{G}(y)
=\{z\}$ or $d_{G}(y)=d_{T_{\{x_2, x_4\}}}(y)\in\{1,2\}$, both contradict $U^4_2\not\subset G$.
This implies that $T_0=\emptyset$.
Again by $U^4_2\not\subset G$, we get every vertex of $T_{x_2}\cup T_{x_4}$ is non-adjacent to $z$. 
Combine this with $T_{x_3}=T_0=\emptyset$ and $v\sim w$, one has $G\subset
H_2$.

$\bullet$ $|T_{\{x_2,x_4\}}|= 1$. We let $T_{\{x_2,x_4\}}= \{w\}$.
We first assert that $T_0=\emptyset$.
Suppose that there exists a vertex, say $v^*$, in $T_0$. In view of Lemma~\ref{c2}(v)-(vi), we have either $d_G(v^*)=d_{T_{\{x_2,x_4\}}}(v^*)=1$ or $d_G(v^*)=d_{T_{\{x_1,x_3\}}}(v^*)=1$, then $U^4_2\subset G$, a contradiction. So we do indeed have $T_0=\emptyset$. 

Since $U^4_2\not\subset G$, by a similar discussion as the proof of Claim~\ref{c9}(i), we immediately get the following two facts.
\begin{fact}\label{f1}
If there exists a number $s\in\{1,3\}$, such that $d_{T_{x_s}}(w)=1$, then
either $T_{x_t}=\emptyset$ or $d_{T_{x_t}}(w)=1$ for $t\in\{1,3\}\setminus \{s\}$.
\end{fact}
\begin{fact}\label{f2}
If there exists a number $p\in\{2,4\}$, such that $d_{T_{x_p}}(z)=1$, then
either $T_{x_q}=\emptyset$ or $d_{T_{x_q}}(z)=1$ for $q\in\{2,4\}\setminus \{p\}$. 
\end{fact}

By $U^4_2\not\subset G$ and Lemma~\ref{c2}(i)-(ii), we immediately get the following two facts.
\begin{fact}\label{f3}
If $d_{T_{x_s}}(w)=1$ for some $s\in\{1,3\}$, then $N_{T_{x_2}}(z)=N_{T_{x_4}}(z)=\emptyset$. 
\end{fact}
\begin{fact}\label{f4}
If $d_{T_{x_p}}(z)=1$ for some $p\in\{2,4\}$, then $N_{T_{x_1}}(w)=N_{T_{x_3}}(w)=\emptyset$. 
\end{fact}
Combine $T_0=\emptyset$ with Facts~\ref{f1}-\ref{f4}, we have either $G\subset
H_{12}$ or $G\subset
H_{13}$.

We now consider the case $|T_{\{x_1,x_3\}}|> |T_{\{x_2,x_4\}}|$. Then $|T_{\{x_1,x_3\}}|=2$ and $|T_{\{x_2,x_4\}}|=1$.
Since $U^4_2\not\subset G$ and Lemma~\ref{c2}(v)-(vi), by a similar discussion as before, we get $T_0=\emptyset$.
Let $T_{\{x_2,x_4\}}= \{w\}$ and let $T_{\{x_1,x_3\}}= \{z_1,z_2\}$.

By $U^4_2\not\subset G$ and Lemma~\ref{c2}(i)-(ii), we get the following fact.
\begin{fact}\label{f5}
$N_{T_{x_1}}(w)=N_{T_{x_3}}(w)=\emptyset$.
\end{fact}

By $Z_3\not\subset G$ and Fact~\ref{f0}, we immediately arrive at the following fact.
\begin{fact}\label{f6}
If there exists a vertex in $T_{x_2}\cup T_{x_4}$, then it is adjacent to exactly one vertex of $T_{\{x_1,x_3\}}$.
\end{fact}

Recall $|T_{x_2}|\geq |T_{x_4}|$. We assert that $|T_{x_2}|=1$ implies $T_{x_4}=\emptyset$. Suppose to the contrary that $|T_{x_2}|=|T_{x_4}|=1$, where $T_{x_2}=\{v\}$ and $T_{x_4}=\{u\}$.
By Fact~\ref{f6}, one has $d_{T_{\{x_1,x_3\}}}(v)=d_{T_{\{x_1,x_3\}}}(u)=1$.
If $N_{T_{\{x_1,x_3\}}}(v)=N_{T_{\{x_1,x_3\}}}(u)$, then $G$ contains $K_5$ as a minor, which is a contradiction.
Thus $N_{T_{\{x_1,x_3\}}}(v)\neq N_{T_{\{x_1,x_3\}}}(u)$, then $U^4_2\subset G$, a contradiction.
Consequently, $|T_{x_2}|=1$ implies $T_{x_4}=\emptyset$. Combine this with $T_0=\emptyset$ and Facts~\ref{f5}-\ref{f6}, one has
$G\subset
H_2$.
This completes the proof of Claim~\ref{c11}.
\end{proof}

{\bf Subcase 3.2.} There exists some $i\in\{1,2,3,4\}$ such that $E(T_{x_i}, T_{x_{i+1}})\neq \emptyset$.

Then there exist vertices $u,v$ such that $u\in T_{x_i}$, $v\in T_{x_{i+1}}$ and $u\sim v$.
By Claim~\ref{c6}(i), we have $|T_2|=1$, then let $T_2=\{w\}$. We may assume that $T_2=T_{\{x_i,x_{i+2}\}}$ (as the case of $T_2=T_{\{x_{i+1},x_{i+3}\}}$ can be prove by symmetry).

By $U^4_2\not\subset G$ and Lemma~\ref{c2}(i),
we immediately get the following fact.
\begin{fact}\label{f7}
$w\sim v$, $w\not\sim u$.
\end{fact}
By $U^4_2\not\subset G$ and $G$ is $K_{3,3}$-minor free, we arrive at Fact~\ref{f8}.
\begin{fact}\label{f8}
$T_{x_{i+3}}=\emptyset$. 
\end{fact}
We distinguish two cases depending on $T_{x_{i+2}}$, then obtain Claims~\ref{c12} and \ref{c13}.

\begin{claim}\label{c12}
If $T_{x_{i+2}}\neq \emptyset$,
then $G\subset
H_{12}$.
\end{claim}

\begin{proof}[\bf Proof of Claim~\ref{c12}]
Since $T_{x_{i+2}}\neq \emptyset$, by the fact that  $G$ is $K_{3,3}$-minor free and $U^4_2\not\subset G$, we have the following fact.
\begin{fact}\label{f9}
There exists exactly one vertex in $T_{x_{i+2}}$ (in what follows, we let $T_{x_{i+2}}=\{z\}$). Moreover, $z\sim v$ and $z\not\sim u$. 
\end{fact}
As $G$ is $K_{3,3}$-minor free, combining Fact~\ref{f9} with Lemma~\ref{c2}(vii) gives the following fact.
\begin{fact}\label{f10}
$T_{x_i}=\{u\}$. 
\end{fact}
Recall that $G$ is $K_5$-minor free and $Z_2, U^4_2\not\subset G$. Using Facts~\ref{f7}-\ref{f10}, Lemma~\ref{c2}(i) and (vii), we arrive at the following fact.
\begin{fact}\label{f11}
$|T_{x_{i+1}}|\leq 2$, where the equality case implies that the vertex of $T_{x_{i+1}}\setminus\{v\}$ is a pendent vertex. %
\end{fact}

Now, combine Facts~\ref{f7}-\ref{f11} with $T_2=\{w\}$, Lemma~\ref{c2}(iii) and (vi), one gets $G\subset
H_{12}$.
This completes the proof of Claim~\ref{c12}.
\end{proof}
Next, we consider the case of $T_{x_{i+2}}= \emptyset$, then get the following claim.

\begin{claim}\label{c13}
If $T_{x_{i+2}}= \emptyset$,
then $G$ is an induced subgraph of some graph of $\{H_{2},H_{11}\}$.
\end{claim}
\begin{proof}[\bf Proof of Claim~\ref{c13}]
Note $Z_3\not\subset G$. In view of Lemma~\ref{c2}(i) and (vii), we get the following fact.
\begin{fact}\label{f13}
$1\leq|T_{x_i}|\leq 2$. Moreover, if $|T_{x_i}|= 2$, then the vertex of $T_{x_i}\setminus\{u\}$ is non-adjacent to any vertex of $\{u,v,w\}$.
\end{fact}

Recall that $v\in T_{x_{i+1}}$. Together with Lemma~\ref{c2}(i), we get $1\leq |T_{x_{i+1}}|\leq 3$, then we consider the following three cases.

$\bullet$ $|T_{x_{i+1}}|=3$.
Then we say $T_{x_{i+1}}\setminus \{v\}=\{v', v''\}$. From Fact~\ref{f7}, we have $w\sim v$.
Since $Z_3\not\subset G$,
$w$ is adjacent to exactly one vertex of $\{v', v''\}$.
Thus, we may assume that $N_{T_{x_{i+1}}}(w)=\{v,v'\}$.

We assert that $T_{x_i}=\{u\}$. Suppose to the contrary that there exists a vertex, say $z$, in $T_{x_i}\setminus\{u\}$. Fact~\ref{f13} gives us that $z$ is non-adjacent to any vertex of $\{u,v,w\}$. 
Since $G$ is $K_{3,3}$-minor free, one has $z\not\sim v'$.
Then $G[\{x_i,x_{i+1},x_{i+2},u,w, v',z\}]\cong Z_3$, a contradiction.
So we do indeed have
$T_{x_i}=\{u\}$.
Combine this with Fact~\ref{f8}, $T_{x_{i+2}}= \emptyset$, $N_{T_{x_{i+1}}}(w)=\{v,v'\}$, Lemma~\ref{c2}(iii) and (vi), we have $G\subset
W_{5}\cong H_2$.

$\bullet$ $|T_{x_{i+1}}|=2$.
Then we say $T_{x_{i+1}}\setminus \{v\}=\{v^*\}$.

If $w\sim v^*$, by a similar discussion as before, one has $T_{x_i}=\{u\}$.
Combine this with Facts~\ref{f7}-\ref{f8}, $T_{x_{i+2}}= \emptyset$, Lemma~\ref{c2}(iii) and (vi), one has $G\subset
W_{5}\cong H_2$.

If $w\not\sim v^*$, by Fact~\ref{f13}, we have $|T_{x_i}|\leq 2$. If there exists a vertex, say $u'$, in $T_{x_i}\setminus \{u\}$, then $u'$ is non-adjacent to any vertex of $\{u,v,w\}$. By $U^4_2\not\subset G$, one has $u'\not\sim v^*$.
Combine this with Lemma~\ref{c2}(iii) and (vi), we get $G\subset
W_{6}\cong H_{11}$.

$\bullet$ $T_{x_{i+1}}=\{v\}$. In view of Fact~\ref{f13}, we get $|T_{x_i}|\leq 2$. Moreover, if $|T_{x_i}|= 2$, then the vertex of $T_{x_i}\setminus\{u\}$ is non-adjacent to any vertex of $\{u,v,w\}$. Together with Lemma~\ref{c2}(iii) and (vi), one has $G\subset
W_{6}\cong H_{11}$.

This completes the proof of Claim~\ref{c13}.
\end{proof}

{\bf Subcase 3.3.} $E(G[T_1])\neq \emptyset$ and $E(T_{x_j}, T_{x_{j+1}})=\emptyset$ for each $j\in\{1,2,3,4\}$. 

In this subcase, there exists some $i\in \{1,2,3,4\}$ such that $e(T_{x_i}, T_{x_{i+2}})\geq 1$. Note that $Z_6\not\subset G$. Hence, we have the following fact immediately.
\begin{fact}\label{f14}
$E(T_{x_{i+1}}, T_{x_{i+3}})=\emptyset$.
\end{fact}

Recall that $E(T_{x_j}, T_{x_{j+1}})=\emptyset$ for every $j\in\{1,2,3,4\}$. Together with $Z_1\not\subset G$ and Lemma~\ref{c2}(i), we immediately get the following fact.
\begin{fact}\label{f15}
$|T_{x_{i+1}}|\leq 1$ and $|T_{x_{i+3}}|\leq 1$.
\end{fact}
By Claim~\ref{c6}(ii), one has $T_2=T_{\{x_i, x_{i+2}\}}$. Since $T_2\neq \emptyset$, there exists a vertex, say $v$, in $T_{\{x_i, x_{i+2}\}}$.
Now by $Z_2\not\subset G$, we immediately get the following fact.
\begin{fact}\label{f16}
Every vertex of $T_{x_i}$ (resp., $T_{x_{i+2}}$) is adjacent to at most one vertex of $T_{x_{i+2}}$(resp., $T_{x_i}$).
\end{fact}

Combine Facts~\ref{f15}-\ref{f16} with Lemma~\ref{c2}(i) and (vii), we get the following fact.
\begin{fact}\label{f17}
$G[T_1\cup \{x_1,x_2,x_3,x_4\}]\subset W_7$.
\end{fact}

We proceed with the following claim.
\begin{claim}\label{c14}
The following statements hold:
\begin{wst}
\item[{\rm (i)}] Every vertex of $T_{x_{i+1}}\cup T_{x_{i+3}}$ is non-adjacent to any vertex of $T_{\{x_i,x_{i+2}\}};$
\item[{\rm (ii)}] Every vertex of $T_0$ is adjacent to exactly one vertex of $T_{\{x_i,x_{i+2}\}}$.
\end{wst}
\end{claim}

\begin{proof}[\bf Proof of Claim~\ref{c14}]
(i) 
Suppose that there exists a vertex $v\in T_{x_{i+1}}\cup T_{x_{i+3}}$, such that $v$ is adjacent to some vertex of $T_{\{x_i,x_{i+2}\}}$, then $U^4_2\subset G$, a contradiction.

(ii) 
By Claim~\ref{c6}(ii) and Lemma~\ref{c2}(v)-(vi), we have $d_G(v)=d_{T_{\{x_i,x_{i+2}\}}}(v)\in \{1,2\}$ for each $v\in T_0.$
Suppose to the contrary that there exists a vertex $v'\in T_0$ such that $d_{T_{\{x_i,x_{i+2}\}}}(v')=2$, by $e(T_{x_i}, T_{x_{i+2}})\geq 1$, we have $U^4_2\subset G$, a contradiction.
This completes the proof of Claim~\ref{c14}.
\end{proof}

Now combine Facts~\ref{f14}-\ref{f17}, Lemma~\ref{c2}(iii), Claim~\ref{c6}(ii) and Claim~\ref{c14}, we have $G\subset \mathcal{F}_1$.

{\bf Case 4.} $T_1=\emptyset$ and $T_2=\emptyset$. 
Then $T_0=\emptyset$, and so $G\cong C_4$, as desired.

This completes the proof of Lemma~\ref{lemg4}.
\end{proof}

We are ready to prove our main result.

\begin{proof}[\bf Proof of Theorem~\ref{Th21}]
{\em Necessity.}
Together with Lemmas~\ref{lemT}, \ref{lemg6}-\ref{lemg4}, our desired result follows immediately.

{\em Sufficiency.} By a direct calculation, we have $\lambda_2(H_i)\leq 1$ for $i\in\{1,2,\ldots,13\}$. Combining with Lemma~\ref{lemF1} gives us that if $G$ is an induced subgraph of a graph of $\mathcal{F}_1\cup\{H_1,H_2,\ldots, H_{13}\}$, then $\lambda_2(G)\leq 1$. The proof is completed.
\end{proof}
\section{\normalsize Concluding remark}\label{s4}
In this paper, we characterize all connected triangle-free planar graphs whose second largest eigenvalue does not exceed $1$, which partially solves Problem~\ref{Pro1}. 

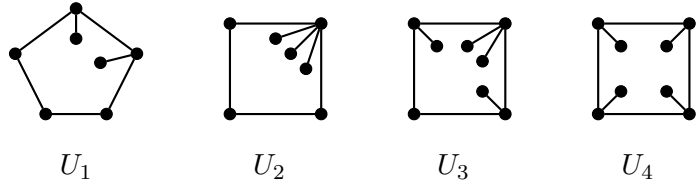
\begin{figure}[!ht]
\centering
  \begin{tikzpicture}[scale = 0.4]
  \tikzstyle{vertex}=[circle,fill=black,minimum size=0.38em,inner sep=0pt]
  \node[vertex] (v_0) at (1,0){};
  \node[vertex] (v_1) at (3,0){};
  \node[vertex] (v_2) at (4,2){};
  \node[vertex] (v_3) at (2,3.5){};
  \node[vertex] (v_4) at (0,2){};
  \node[vertex] (v_5) at (2,2.5){};
  \node[vertex] (v_6) at (2.8,1.7){};
  \draw[fill](1,0)circle(.19); 
  \draw[fill](3,0)circle(.19); 
  \draw[fill](4,2)circle(.19); 
  \draw[fill](2,3.5)circle(.19); 
  \draw[fill](0,2)circle(.19); 
  \draw[fill](2,2.5)circle(.19); 
  \draw[fill](2.8,1.7)circle(.19); 
  \draw[thick] (v_0) -- (v_1)--(v_2)--(v_3)--(v_4)--(v_0);
  \draw[thick] (v_3) -- (v_5);
  \draw[thick] (v_2) -- (v_6);
 
 \node[below] at (2,-1){$U_1$};
  \end{tikzpicture}
 \hspace{0.8cm}
  \begin{tikzpicture}[scale = 0.4]
  \tikzstyle{vertex}=[circle,fill=black,minimum size=0.38em,inner sep=0pt]
  \node[vertex] (v_0) at (0,0){};
  \node[vertex] (v_1) at (3,0){};
  \node[vertex] (v_2) at (3,3){};
  \node[vertex] (v_3) at (0,3){};
  \node[vertex] (v_4) at (1.5,2.5){};
  \node[vertex] (v_5) at (2,2){};
  \node[vertex] (v_6) at (2.5,1.5){};
  \draw[fill](0,0)circle(.19); 
  \draw[fill](3,0)circle(.19); 
  \draw[fill](3,3)circle(.19); 
  \draw[fill](0,3)circle(.19); 
  \draw[fill](1.5,2.5)circle(.19); 
  \draw[fill](2,2)circle(.19); 
  \draw[fill](2.5,1.5)circle(.19); 
  \draw[thick] (v_0) -- (v_1)--(v_2)--(v_3)--(v_0);
  \draw[thick] (v_2) -- (v_4);
  \draw[thick] (v_2) -- (v_5);
  \draw[thick] (v_2)--(v_6);
 \node[below] at (1.3,-1){$U_2$};
  \end{tikzpicture}
 \hspace{0.8cm}
  \begin{tikzpicture}[scale = 0.4]
  \tikzstyle{vertex}=[circle,fill=black,minimum size=0.38em,inner sep=0pt]
  \node[vertex] (v_0) at (0,0){};
  \node[vertex] (v_1) at (3,0){};
  \node[vertex] (v_2) at (3,3){};
  \node[vertex] (v_3) at (0,3){};
  \node[vertex] (v_4) at (1.75,2.25){};
  \node[vertex] (v_5) at (0.75,2.25){};
  \node[vertex] (v_6) at (2.25,1.75){};
  \node[vertex] (v_7) at (2.25,0.75){};
  \draw[fill](0,0)circle(.19); 
  \draw[fill](3,0)circle(.19); 
  \draw[fill](3,3)circle(.19); 
  \draw[fill](0,3)circle(.19); 
  \draw[fill](1.75,2.25)circle(.19); 
  \draw[fill](0.75,2.25)circle(.19); 
  \draw[fill](2.25,1.75)circle(.19); 
  \draw[fill](2.25,0.75)circle(.19); 
  \draw[thick] (v_0) -- (v_1)--(v_2)--(v_3)--(v_0);
  \draw[thick] (v_2) -- (v_4);
  \draw[thick] (v_3) -- (v_5);
  \draw[thick] (v_2)--(v_6);
  \draw[thick] (v_1)--(v_7);
 \node[below] at (1.3,-1){$U_3$};
  \end{tikzpicture}
 \hspace{0.8cm}
  \begin{tikzpicture}[scale = 0.4]
  \tikzstyle{vertex}=[circle,fill=black,minimum size=0.38em,inner sep=0pt]
  \node[vertex] (v_0) at (0,0){};
  \node[vertex] (v_1) at (3,0){};
  \node[vertex] (v_2) at (3,3){};
  \node[vertex] (v_3) at (0,3){};
  \node[vertex] (v_4) at (2.25,2.25){};
  \node[vertex] (v_5) at (0.75,2.25){};
  \node[vertex] (v_6) at (0.75,0.75){};
  \node[vertex] (v_7) at (2.25,0.75){};
  \draw[fill](0,0)circle(.19); 
  \draw[fill](3,0)circle(.19); 
  \draw[fill](3,3)circle(.19); 
  \draw[fill](0,3)circle(.19); 
  \draw[fill](2.25,2.25)circle(.19); 
  \draw[fill](0.75,2.25)circle(.19); 
  \draw[fill](0.75,0.75)circle(.19); 
  \draw[fill](2.25,0.75)circle(.19); 
  \draw[thick] (v_0) -- (v_1)--(v_2)--(v_3)--(v_0);
  \draw[thick] (v_2) -- (v_4);
  \draw[thick] (v_3) -- (v_5);
  \draw[thick] (v_0)--(v_6);
  \draw[thick] (v_1)--(v_7);
 \node[below] at (1.3,-1){$U_4$};
  \end{tikzpicture}
\caption{Graphs $U_1,\ldots,U_4$.}\label{fG}
\end{figure}

In fact, our obtained result has closely relationship with some known results. We first recall some definitions. Let $G$ be an $n$-vertex connected graph. We call $G$ a \textit{unicyclic graph} if the size of $G$ is $n$; whereas we call $G$ a \textit{bicyclic graph} if the size of $G$ is of $n+1$. 
Let $K^3_{s_1;s_2;s_3}$ be the graph obtained by attaching $s_i$ pendant edges to each vertex of $K_3,$ where $s_i\geqslant 0$.
Let $\mathcal{U}$ be the class of graphs each of which is obtained by attaching $t$ pendant paths of length 2 to a vertex of $K_3$, where $t\geqslant 0.$
Let $U_1,\ldots,U_4$ be the unicyclic graphs as depicted in Figure~\ref{fG}. 
In what follows we will illustrate the relation between Theorem~\ref{Th21} with some known results.
\begin{thm}[Xu~\cite{XuGH1}]
For a unicyclic graph $G$, $\lambda_2(G)\le 1$ holds if and only if $G$ is an induced unicyclic subgraph of a graph belonging to $\mathcal{U}\cup\{C_6,K^3_{3;1;1},K^3_{4;1;0},K^3_{2;2;2},U_1, \ldots,U_4\}.$
\end{thm}

Note that $C_6\subset \mathcal{F}_1$, $U_1\subset \mathcal{F}_1$, $U_2\subset H_6$, $U_3\subset H_6$ and $U_4\subset H_8$. Thus Theorem~\ref{Th21} deduces all triangle-free unicyclic graphs with $\lambda_2\le 1.$ 

Let $G_1,\ldots,G_{14}$ be the bicyclic graphs illustrated in \cite[Fig.~1]{GuoSG}.
\begin{thm}[Guo~\cite{GuoSG}]
For a bicyclic graph $G$, $\lambda_2(G)\le 1$ holds if and only if $G$ is an induced bicyclic subgraph of a graph belonging to $\{G_1, \ldots,G_{14}\}.$
\end{thm}

Note that $G_i\subset H_5$ for $i\in\{9,10\}$, $G_j\subset \mathcal{F}_1$ for $j\in\{11,\ldots,14\}$ and $G_k$ contains triangle for $k\in\{1,\ldots,8\}.$ Therefore, our result may derives all triangle-free bicyclic graphs with $\lambda_2\le 1.$ 

A graph consisting of $t$ internal disjoint $(u,v)$-paths is called {\em generalized $\theta$-graph}, where $t\ge 2.$ 
Let $\Gamma_1$ and $\Gamma_4$ be the generalized $\theta$-graphs illustrated in \cite[Fig.~1]{HQX}. 

\begin{thm}[Gao and Huang~\cite{HQX}]
For a generalized $\theta$-graph $G$, $\lambda_2(G)\le 1$ holds if and only if $G$ is an induced generalized $\theta$-subgraph of $\Gamma_1$ or $\Gamma_4$.
\end{thm}

Note that $\Gamma_1\subset \mathcal{F}_1,$ and if $G$ is an induced generalized $\theta$-subgraph of $\Gamma_4$ without triangle, then $G\subset H_5$. Therefore, our result may deduce all triangle-free generalized $\theta$-graphs with $\lambda_2\le 1.$

The following problem is interesting and challenging.
\begin{problem}
    Characterize all connected planar graphs of girth $3$ whose second largest eigenvalue does not exceed $1$.
\end{problem}

\section*{\normalsize Acknowledgement}
The authors would like to thank the referees for their valuable comments which lead to an improvement of the original manuscript.


\end{document}